\documentclass[11pt]{article}
\usepackage {amssymb}
\usepackage {amsmath}
\usepackage {amsfonts}
\newcommand{\OPS}{\mathrm{OPS}}

\def\<{{\langle}}
\def\>{{\rangle}}
\def\a{{\alpha}}
\def\C{{\mathbb{C}}}
\def\l{{\lambda}}

\def\b{{\beta}}

\def\Z{{\Sigma}}

\def\1{{\Sigma_1}}
\def\2{{\Sigma_2}}
\def\R{{\mathbb{R}}}
\def\S{{\mathbb{S}}}
\def\N{{\mathbb{N}}}

\def\g{{\gamma}}
\def\e{{\epsilon}}
\def\tr{\textrm}
\def\v{{\varrho}}
\def\p{{\partial}}

\def \fine{\hfill {\rule {2mm} {2mm}} \vspace{3mm}}

\def\ds{\displaystyle}
\def\mr{\mathrm}
\newtheorem{theorem}{Theorem}[section]
\newtheorem{lemma}[theorem]{Lemma}
\newtheorem{proposition}[theorem]{Proposition}
\newtheorem{definition}{Definition}[section]

\newtheorem{remark}[theorem]{\emph{Remark}}
\begin{document}
\numberwithin{equation}{section}
\title{Hypoellipticity and Higher Order Levi Conditions}
\author{Marco Mughetti}
\date{\small Department of Mathematics, University of Bologna,\\ Piazza di Porta S.\ Donato 5, 40127 Bologna, Italy\\
e-mail: marco.mughetti@unibo.it\\} \maketitle
\begin{abstract} We study the $C^\infty$-hypoellipticity for a class of double characteristic operators with simplectic characteristic manifold, in the case the classical condition of minimal loss of derivatives is violated.
\end{abstract}
\section{Introduction}
The purpose of this paper is to establish necessary and sufficient conditions for the $C^\infty$ regularity of the solutions of the Grushin operators of the
following form
\begin{equation}\label{grushin}
 P=D_{x_1}^2+a(x_1,x'){x_1}^{2h}\Delta_{x'}+{x_1}^{h-1}B(x_1,x',D_{x'}),\qquad x_1\in\R,\;x'=(x_2,...,x_{n})\in\R^{n-1},
\end{equation}
where $D=\frac{1}{i}\partial$, $h$ is a positive integer, $a(x_1,x')$ is a strictly positive smooth function in $\R^n$, $\Delta_{x'}$ is the positive Laplacian in
$\R_{x'}^{n-1}$ and $B(x_1,x',D_{x'})$ is a classical first order (pseudo)differential operator in $x'$ depending on the parameter $x_1$. The first order term of
$P$ is chosen in such a way to behave as the commutators $[D_{x_1},x_1^hD_{x_j}]$ of the vector fields $D_{x_1}$ and $x_1^hD_{x_j}$ ($j=2,...,n$), generating its
principal part. We are interested in studying the effect of the lower order terms in $B(x_1,x',D_{x'})$ on the hypoellipticity of $P$.\\
We point out that the Grushin model \eqref{grushin} represents a typical canonical form of the (pseudo)-differential operators with symplectic characteristic
manifold of codimension $2$.\\ As we shall make clear later on, the $C^\infty$-hypoellipticity of such operators is strictly related to the spectral properties of
an anharmonic oscillator (see \eqref{locgrushin} below), whose spectrum, unlike the case $h=1$, cannot be explicitly computed.\\
%In Section \ref{} we use the perturbation theory of linear operators (see \cite{kato}) to overcome this difficulty.\\
The problem of $C^\infty$-hypoellipticity has been widely studied in literature (see, for instance, \cite{ho} Chapters XXII and XXVII and the references therein).
We recall that a (pseudo)differential operator $Q(x,D_x)$ is $C^\infty$-hypoelliptic (h.e.) if $Q$ preserves the $C^\infty$ singular support, i.e.
$$
\mathrm{sing\,\ supp}\; Qu=\mathrm{sing\,\ supp}\; u,\qquad
\mathrm{for\,\ any }\;u\in{\mathcal{D'}}.
$$
This definition can be made more precise if we introduce the
notion of loss of derivatives; namely, we say that a
pseudodifferential operator $Q$ of order $m$ in an open set
$X\subset\R^n$ is $C^\infty$-hypoelliptic, with loss of
$\gamma\geq 0$ derivatives, if
\begin{equation}\label{defipo}
\textrm{\bf(He)}_\g\quad \forall s\in\R, \forall u\in{\cal D}'(X),
\forall\,\mathrm{open\; set}\; \Omega\,\subset X: Qu\in H^s_{\rm
loc}(\Omega) \Longrightarrow u\in H^{m+s-\gamma}_{\rm
loc}(\Omega).
\end{equation}
A well-known sufficient condition for the h.e. of (\ref{grushin}) is the injectivity in $L^2(\R)$ of the following differential operator in the $t\in\R$ variable
(the so-called localized operator associated with $P$)
\begin{equation}\label{locgrushin}
P_{\varrho}=D_t^2+a(0,x')t^{2h}|\xi'|^2+t^{h-1}b_1(0,x',\xi')
\end{equation}
as $\varrho=(x',\xi')$ varies in $T^\ast\R^{n-1}\setminus 0$,  $b_1(x,\xi')$ being the principal symbol of $B$
(see \cite{gr}, \cite{grushin}, \cite{gilioli2}, \cite{menikoff} and also Thm.1.12 \cite{marco}).
In this case $P$ turns out to be hypoelliptic with loss of $2h/(h+1)$ derivatives, and this is exactly the minimal loss of regularity  we can expect from $P$.\\
This assumption however is far from being necessary; we can, in general, have that $P$ is hypoelliptic although $P_{(x'_0,\xi'_0)}$ is not injective in $L^2(\R)$ at some point $(x'_0,\xi'_0)$. The present paper is  devoted to the study of the $C^\infty$-hypoellipticity of the operator $P$, for which the classical $L^2-$injectivity condition is violated at some point of its characteristic manifold.\\
We point out that our results yield an alternative proof of the $C^\infty$-hypoellipticity of the Kohn operator  introduced in  \cite{kohn} and \cite{bove1}.
Actually we discuss the following extension of the kohn model:
\begin{equation}\label{kintro}
LL^*+\big(f(x_1)L\big)^*f(x_1)L,\qquad L=D_{x_1}+i g(x_1)D_{x_2}
\end{equation}
where $f(x_1), g(x_1)$ are  polynomial  functions, so that $L$ can be regarded as a generalization of the Mizohata operator $M=D_{x_1}-i x_1^h D_{x_2}$ ($h\in\N$).
We shall show that, in that case, the hypoellipticity strictly depends on the presence of common zeroes of the functions $g(x_1), f(x_1)$ and on their order (see
Proposition \ref{hkext}). As pointed out in  \cite{kohn}, the loss of derivatives of \eqref{kintro} can be very large compared to the case of sums of squares of
real vector fields where the loss is always strictly less than $2$ (see Thm. 22.2.1 \cite{ho} and  \cite{RS}). In this connection, in Section \ref{uniquefield}
(Prop. \ref{propunique}) we show that this phenomenon takes place also in the presence of a single  complex vector field, by analyzing the following example:
$$\big(D_{x_1}-i x_1^{k_1}D_{x_2}\big)^\ast \big(D_{x_1}-i x_1^{k_1}D_{x_2}\big)+\big(x_1^{k_2} D_{x_2})^\ast \big(x_1^{k_2} D_{x_2}).$$

If $h=1$, the principal symbol of \eqref{grushin} vanishes exactly to the second order on its characteristic manifold (the so-called transversally elliptic case).
This case has been extensively studied over the years, see for instance \cite{boutet}, \cite{h4}, \cite{grushin}, \cite{h2}, \cite{kw}, \cite{Parenti3}. All the
above mentioned articles concern only transversally elliptic operators. Here we consider the Grushin-type operators \eqref{grushin}, where the transversal ellipticity fails because of $h>1$. We study the $C^\infty$ hypoellipticity of $P$ assuming that its localized operator \eqref{locgrushin} is not injective at some
point $(x_0',\xi_0')$, or, equivalently, there exists an eigenvalue $\l_{j_0}$ of \eqref{locgrushin} such that $\l_{j_0}(x_0',\xi_0')=0$. As we shall see later on
(Theorem \ref{teoprinc}), in this framework the minimal loss of regularity $2h/(h+1)+\delta(h)$ of $P$ strictly depends on the parity of $h$. Furthermore, in the
transversally elliptic case, in \cite{h2} it is shown that the minimal loss is $1+1/2$ and is attained iff attained if and only if
\begin{equation}\label{autopoisson2}
\frac{1}{i}\{\l_{j_0},\overline{\l_{j_0}}\}(x_0',\xi_0')<0,
\end{equation}
where $\{\, ,\, \}$ denote the Poisson brackets in $T^\ast\R^{n-1}$, and $\overline{\l_{j_0}}$ the complex conjugate of $\l_{j_0}$. In the non transversally elliptic case ($h>1$) the situation is much more delicate. Let us consider for instance relation \eqref{autopoisson2}: we shall see that it does not play any role if $h\neq 3$. However, also for the ``special value" $h=3$, \eqref{autopoisson2} is generally no longer a necessary and sufficient condition for the $C^\infty$ hypoellipticity of $P$ (see Section \ref{mod}).\\
The situation is quite different if the coefficients of $P$ only
depend on the tangent variables to the characteristic
manifold $\Sigma=\{x_1=0=\xi_1\}$, say $x'$,
\begin{equation}\label{grushintan}
D_{x_1}^2+a(x'){x_1}^{2h}\Delta_{x'}+{x_1}^{h-1}B(x',D_{x'}),\qquad
x_1\in\R,\;x'\in\R^{n-1}.
\end{equation}
In this particular setting, the minimal loss of regularity is $2h/(h+1)+1/2$ and is attained only for $h$ odd, under the classical hypothesis \eqref{autopoisson2},
$\l_{j_0}$ being the eigenvalue of \eqref{locgrushin} that vanishes at $(x_0',\xi_0')$. If $h$ is an even integer, the minimal loss of derivatives must necessarily be larger than $2h/(h+1)+1/2$ (see Section \ref{mod}). We stress the fact that $\l_{j_0}$ is not explicitly known if $h>1$; nevertheless, we shall show that condition \eqref{autopoisson2} can be directly deduced from \eqref{grushintan} by using some results of Perturbation theory.

% The aim of this paper is to study the $C^\infty$-hypoellipticity of the Grushin-type operators \eqref{grushin}, for which the transversal elliptic condition fails (that is $h>1$) and its localized operator \eqref{locgrushin} is not injective at some point $(x_0',\xi_0')$.
%We shall show that the minimal loss of regularity $2h/(h+1)+\delta$ of $P$ strictly depends on the parity of $h$, and the relation \eqref{autopoisson} is no longer a necessary and sufficient condition.
%unlike the case $h=1$,  $\delta$ is generally smaller than $1/2$.
The paper is organized as follows: in Section \ref{machinery} we introduce the machinery required to study the hypoellipticity of the operator \eqref{grushin},
microlocally near a degenerate point; in Section \ref{ko} we apply our results to some classes of examples; in Section \ref{mod} we analyse the h.e. of
Grushin operator of type (\ref{grushintan}). Finally, in Appendix I we discuss the hypoellipticity properties of a class of one-dimensional polyhomogeneous symbols (see Def.
18.1.5 \cite{ho}) naturally associated with the operator \eqref{grushin} via Theorem \ref{teoprinc}; in Appendix II we develop a suitable anisotropic
variant of the Boutet de Monvel pseudodifferential calculus introduced in \cite{boutet} and in \cite{h2}.

%\par\medskip{\bf Acknowledgments.}

\section{The general case}\label{machinery}
Consider the operator $P$ in \eqref{grushin} and its localized operator $P_\varrho$ in \eqref{locgrushin}. The characteristic manifold $\Z=\{(x,\xi)\in T^\ast\R^n\setminus 0\, |\,x_1=\xi_1=0\} $ of $P$ can be trivially identified with $T^\ast\R^{n-1}\setminus 0$.\\
It is well-known that the spectrum of $P_\varrho$ is discrete,
i.e. $\textrm{Spec}(P_\varrho)$ is made of simple isolated
eigenvalues $\l_j(\varrho), j\in\N,$ of finite multiplicities,
diverging to $+\infty$ (see, for instance,Theorem 3.1 and
Proposition 3.3 in Chapter 2 \cite{bs}):
\begin{equation}\label{ord}
\l_0(\varrho)<\l_1(\varrho)<...<\l_j(\varrho)<\l_{j+1}(\varrho)<...,\qquad
j\in\N.
\end{equation}
In the classical study of the $C^\infty$-hypoellipticity of $P$ (see \cite{boutet,h4, gr, grushin, marco}) the key point concerns the assumption $\l_j(\v)\neq 0$
for every $j=0,1,2,...$; as a consequence, $P_\v$ turns out to have a suitable left inverse, by means of which one can construct a left parametrix of $P$ via the
Boutet de Monvel's pseudodifferential calculus (see, for instance, Appendix II); this yields the $C^\infty$-hypoellipticity of $P$ with minimal loss of $2h/(h+1)$
derivatives. On the other hand, if there exists a point $\varrho=\varrho_0:=(x'_0,\xi'_0)$ and an eigenvalues $\l_{j_0}$ such that
\begin{equation}\label{avr}
\l_{j_0}(x'_0,\xi'_0)=0,
\end{equation}
$P$ can be $C^\infty$-hypoelliptic only with loss of derivatives larger than $2h/(h+1)$.\\
In order to introduce our results,  let us briefly discuss the wellknown case $h=1$, for which the localized operator $P_{(x',\xi')}$ in \eqref{locgrushin} is a
harmonic oscillator and its eigenvalues $\{\l_j(x',\xi')\}_{j\in\N}$ can be explicitly computed. Due to \cite{gr}, \cite{grushin}, \cite{boutet}, $P$ is h.e. with
loss of one derivate iff $P_{(x',\xi')}$  is $L^2(\R)-$injective for every $(x',\xi')\in T^\ast\R^{n-1}\setminus 0$, which means that, for any $(x',\xi')\in
T^\ast\R^{n-1}\setminus 0$
\begin{equation}\label{violate}
0\neq\l_j(x',\xi')=(2j+1)\sqrt{a(0,x',\xi')}+b_1(0,x',\xi'),\qquad \forall j\in\N.
\end{equation}
In \cite{h2} Helffer assumed that condition (\ref{violate}) fails at some point $(x_0',\xi_0')$, yielding \eqref{avr} or, equivalently, that
$b_1(0,x_0',\xi_0')=-(2j_0+1)\sqrt{a(0,x',\xi')}$. In this framework, $P$ can be hypoelliptic only with loss $r\geq 1+1/2=3/2$ derivatives (see \cite{sj});
furthermore, in Theorem 1.2 \cite{h2} it is shown that the bound $3/2$ is attained if and only if
\begin{equation}\label{autopoisson}
\frac{1}{i}\{\l_{j_0},\overline{\l_{j_0}}\}(x_0',\xi_0')<0.
\end{equation}
Parenti and Parmeggiani \cite{Parenti3} have pushed forward this analysis, providing necessary and sufficient conditions for the h.e. (with a big loss of
regularity) of a large class of transversally elliptic Grushin operators.
The same problem was also studied by Kwon \cite{kw} by using the Treves concatenations method (see \cite{tre}).\\
Suppose now that $h>1$, the first difficulty we face up  concerns the fact that the localized operator \eqref{locgrushin} is an anharmonic-type oscillator and hence
its eigenvalues $\{\l_j(x',\xi')\}_{j\in\N}$ cannot be explicitly computed. Nevertheless, in Theorem \cite{marco2} it is proved that \eqref{avr} amounts to saying
that
\begin{equation}\label{avreven}
|b_1(0,x'_0,\xi'_0)|=\sqrt{a(0,x'_0)}(h+1)(2j_0+1)|\xi'_0|\quad\textrm{if $h$ is even};
\end{equation}
\begin{equation}\label{avrodd}
b_1(0,x'_0,\xi'_0)=(-1)^{j_0}\sqrt{a(0,x'_0)}|\xi'_0|-\sqrt{a(0,x'_0)}(h+1)(j_0+\theta(j_0))|\xi'_0|\quad\textrm{if $h$ is odd},
\end{equation}
where $\theta(j):=1$ if $j$ is even and $\theta(j):=0$ if $j$ is odd.\\
Our aim is to show that the h.e. of $P$ is equivalent to the h.e. of an operator $L(y,D)\in\textrm{OP}S^{m-2h/(h+1)}(\R^{n-1})$
in fewer variables, easier to be analysed.\\
In our setting, $P_{\v_0}$ is clearly not invertible, hence we
cannot use the classical approach described above. However, we can
exploit the following linear algebra remark: assume that the $ n
\times n $ matrix $ A $ has zero in its spectrum with multiplicity
one. Then of course $ A $ is not invertible, but, denoting by $
e_{0} $ the zero eigenvector of $ A $, the matrix (in block form)
$$
\begin{bmatrix}
A & e_{0} \\
{}^{t}e_{0} & 0
\end{bmatrix}
$$
is invertible as a $ (n+1) \times (n+1) $ matrix in $ \C^{n+1} $.
Here $ {}^{t}e_{0} $ denotes the row vector $ e_{0} $. This
strategy goes back to Grushin \cite{gr} and Sj\"ostrand \cite{sj2}
and then exploited by Helffer \cite{h2}, by Grigis-Rothschild
\cite{grigis}, by Parenti-Parmeggiani
\cite{Parenti3} and, recently, by Bove-Mughetti-Tartakoff \cite{bmt1}, \cite{bmt2} (see also \cite{sj3} for a complete survey).
Roughly speaking, the idea is to replace the operator $P$ by a
suitably chosen square system of operators
\begin{equation}\label{sist}
{\cal A}= \left(\begin{array}{cc}
\ds P & H^-\\
H^+ & 0
\end{array}\right)
\end{equation}
in such a way that its ``localized operator"
\begin{equation}\label{locsis}
{\cal A_\varrho}= \left(\begin{array}{cc}
\ds P_\varrho & h_{\phi_2}^-(\v)\\
h_{\phi_1}^+(\v) & 0
\end{array}\right)
\end{equation}
turns out to be injective (actually invertible) for any $\varrho\in\Z$, although this is clearly false for the localized operator $P_{\varrho}$ in (\ref{locgrushin}) at $\v=\v_0$ because of \eqref{avr}.\\
As a consequence of the calculus we shall develop, it follows that
$\cal A$ admits a two-sided parametrix
\begin{equation}\label{param}
{\cal E}= \left(\begin{array}{cc}
\ds E & K^-\\
K^+ & -L
\end{array}\right)
\end{equation}
where $L=L(y,D_y)$ is exactly the operator we were seeking for. We emphasize that $\cal A_\varrho$ is invertible at every $\varrho\in\Z$, and
its inverse can be explicitly computed without using Neumann series; this fact is crucial if you need to know  the complete symbol of $L$ and
is the reason we cannot proceed as in \cite{h2}. Here we prefer to follow the approach used in \cite{Parenti3}.\\
The key points in the above program are the choice of the operators $H^\pm$ in (\ref{sist}) and the construction of the pseudodifferential calculus (i.e., the
classes of symbols and the related composition rules) based on them and on the operator $P$ in \eqref{grushin}. This task is very technical and here we prefer to
develop only the crucial points of the machinery we need, putting in evidence  the required adjustments and referring the reader to \cite{Parenti3} and to \cite{bmt2} for further details.\\
Our starting point is the construction of the operators $h_{\phi_1}^+(\v), h_{\phi_2}^-(\v)$  (and afterwards of the related operators $H^\pm$) in such a way that
the localized system ${\cal A_\varrho}$ is invertible at any $\v\in\Z$.\\
Let us consider the localized operator
$P_\v:L^2(\R)\longrightarrow L^2(\R)$ as an unbounded operator
with the ``natural'' domain
$$
B_h^2(\R)=\{f\in{\mathcal S}'(\R):
\|f\|_{B^2_h}=\Big(\sum_{\a/h+\b\leq 2}\|t^\a D_t^\b
f\|^2_{L^2(\R)}\Big)^\frac{1}{2}<+\infty\}.$$
The null eigenfunctions of $P_{\v_0}$ plays a fundamental role in the construction of the operators $h_{\phi_1}^+(\v), h_{\phi_2}^-(\v)$ in \eqref{locsis};
roughly speaking, they are obtained by using  suitable microlocal extensions $\phi_1(\v,.), \phi_2(\v,.)$, near $\v_0$, of such null eigenfunctions.\\
Precisely, since the principal symbol $p_2$ of $P$ is assumed to be real, from \eqref{avreven}, \eqref{avrodd} it follows that all the coefficients of $P_{\v_0}$ in
(\ref{locgrushin}) are real; hence, $P_{\varrho}$ turns out to be a self-adjoint operator at $\varrho=\varrho_0$. Moreover, every eigenvalue in \eqref{ord} is
simple (see Proposition 3.3 Chap.2 \cite{bs}) so that the kernel of $P_{\v_0}=P_{\v_0}^\ast$ is one-dimensional, i.e. there exists a function $0\neq\phi(\v_0; .
)\in{\mathcal S}(\R)$ for which $\textrm{Ker}\ P_{\v_0}=\langle \phi(\v_0; . )\rangle$. The main point is now the way we extend $\phi(\v_0; . )$ microlocally near
$\v_0$ in order to get the functions $\phi_1(\v,.), \phi_2(\v,.)$.
In \cite{h2} $\phi_1(\v,.), \phi_2(\v,.)$ are chosen equal to the $L^2(\R)$-normalized  eigenfunction associated with the eigenvalue $\l_{j_0}(\v)$ (see \eqref{avr}). Unfortunately, it seems to us that this choice does not allow to get an explicit inverse of ${\cal A_\varrho}$ near $\v_0$, and we thus follows a slight different approach.\\
Let us start off by considering the square operators $P_{\v}^\ast
P_{\v}$ and $P_{\v} P_{\v}^\ast$; they are $h-$globally elliptic,
self-adjoint, non negative differential operators (see \cite{h3}
and Section 1.5 \cite{marco} for the $h$-anisotropic version of the
pseudodifferential calculus), with discrete spectrum contained in
$[0,+\infty[$. Denote by $\mu_1(\v), \mu_2(\v)$ the smallest
eigenvalue of $P_{\v}^\ast P_{\v},\, P_{\v} P_{\v}^\ast$,
respectively, and consider the related eigenspaces $V_1(\v):=
\textrm{Ker}\ \big(P_{\v}^\ast P_{\v}-\mu_1(\v)I\big)$ and
$V_2(\v):= \textrm{Ker}\ \big(P_{\v} P_{\v}^\ast-\mu_2(\v)I\big)$.
We have the following result.
\begin{lemma}\label{quad}
There exist a conic neighborhood $U\subset\Z\cong
T^\ast\R^{n-1}\setminus 0$ of $\v_0=(x'_0,\xi'_0)$ and two
functions $\phi_1(\v;t), \phi_2(\v;t)\in C^\infty\big(U,{\mathcal
S}(\R)\big)$ such that, for any $\v=(x',\xi')\in U$
\begin{enumerate}
%\item $U\ni\v=(x',\xi')\rightarrow \l_i(\v)$ are smooth maps for $i=1,2$;
\item $V_1(\v)=\langle \phi_1(\v; \, . \,)\rangle$,
$V_2(\v)=\langle \phi_2(\v; \, . \,)\rangle;$ \item
$\phi_i(x',s\xi',s^{-\frac{1}{h+1}}t)=s^{\frac{1}{2(h+1)}}
\phi_i(x',\xi',t),\,s>0,\;i=1,2;$ \item if $h$ is an odd integer,
$\phi_1(\rho;t)$ and $\phi_2(\rho;t)$ are even or odd functions in
the $t$ variable, with the same parity; namely, for $i=1,2$,
    $$\phi_i(\rho;-t)=\phi_i(\rho;t),\;\;\forall \rho\in U, \quad{ or}\quad\phi_i(\rho;-t)=-\phi_i(\rho;t), \;\;\forall \rho\in U;$$

\item $\ds\int |\phi_i(x',\xi',t)|^2 dt=1,\;i=1,2.$
\end{enumerate}
\end{lemma}
{\bf Proof.}
We shall prove the two statements only for the eigenvalue $\l_1$ and the space $V_1$; the same arguments apply to $\l_2$ and $V_2$.\\
To begin with, in view of the homogeneity properties of $P_\v$, it is enough to prove the above statements in a neighborhood of $\v_0$ in
$\mathbb{S}^*\Z=\{\v=(x',\xi')\in\Z, |\xi'|=1\}$.\\
To see this, note that the symbol of the localized operator $P_{(x',\xi')}$ (see \eqref{locgrushin}) satisfies the following property
$$p_{(x',s\xi')}(s^{-\frac{1}{h+1}}t,s^{\frac{1}{h+1}}\tau)=s^{\frac{2}{h+1}}
p_{(x',\xi')}(t,\tau),\quad\textrm{for any}\, s>0.$$ Upon considering the unitary operator
$$
M_s : L^2(\R) \longrightarrow L^2(\R),\quad (M_sf)(t) =
s^{\frac{1}{2(h+1)}} f(s^{\frac{1}{h+1}}t),\quad s > 0,
$$
it is easily seen that
$$P_{(x',s\xi')}(M_sf)=s^{\frac{2}{h+1}}M_s(P_{(x',\xi')}f),$$
and a similar computation for the adjoint operator
$P_\varrho^\ast$ yields
$$P_{(x',s\xi')}^\ast(M_sf)=s^{\frac{2}{h+1}}M_s(P^\ast_{(x',\xi')}f).$$
As a consequence we obtain
\begin{equation}\label{homp}
P_{(x',s\xi')}^\ast
P_{(x',s\xi')}(M_sf)=s^{\frac{4}{h+1}}M_s(P_{(x',\xi')}^\ast
P_{(x',\xi')}f).
\end{equation}
Since
$$
\mu_1(x',\xi')=\min_{f\in{\mathcal S}(\R), \|f\|_{L^2}=1}\;\langle
P_{(x',\xi')}^\ast P_{(x',\xi')}f,f\rangle,
$$
from \eqref{homp} it follows that
\begin{equation}\label{homa}
\mu_1(x',s\xi')=s^{\frac{4}{h+1}}\mu_1(x',\xi')
\end{equation}
whence
\begin{equation}\label{autohom}
\textrm{Ker}\big(P_{(x',\xi')}^\ast
P_{(x',\xi')}-\mu_1(x',\xi')\big)=M_{s^{-1}}\Big(\textrm{Ker}\big(P_{(x',s\xi')}^\ast
P_{(x',s\xi')}-\mu_1(x',s\xi')\big)\Big),\, s>0.
\end{equation}
We thus have $V_1(x',\xi')=M_{|\xi'|}\big(V_1(x',\xi'/|\xi'|)\big)$.\\
Therefore it suffices to prove that the statements of Lemma \ref{quad} hold in a neighborhood $\tilde{U}\subset \mathbb{S}^*\Z$ of $\v_0\in\mathbb{S}^*\Z$.\\
The lowest eigenvalue $\mu_1(\v)$ is, a priori, only a  continuous function on $\v\in\Z$. However, since the operators $P_{\v}^\ast P_{\v}$ has discrete spectrum,
there exist a connected neighborhood $\tilde{U}\subset \mathbb{S}^*\Z$ of $\v_0$ and a constant $0<\delta\in\R$ such that $[0,\delta]\cap\textrm{Spec}\,
(P_{\v}^\ast P_{\v})=\{\mu_1(\v)\}$ for all $\v\in\tilde{U}$. By standard spectral theory (see \cite{kato} or \cite{RS} Chapter XII) the orthogonal projector
$\pi_1(\v)$ onto $V_1(\v):= \textrm{Ker}\ \big(P_{\v}^\ast P_{\v}-\mu_1(\v)I\big)$ is given by
\begin{equation}\label{pro1}
\pi_1(\v)=\frac{1}{2\pi i}\oint_{|\zeta|=\delta}(\zeta-P_{\v}^\ast
P_{\v})^{-1}d\zeta, \qquad \forall\v\in\tilde{U},
\end{equation}
and, similarly, by possibly shrinking the neighborhood
$\tilde{U}$, for the orthogonal projector $\pi_2(\v)$ onto
$V_2(\v):= \textrm{Ker}\ \big(P_{\v} P_{\v}^\ast-\mu_2(\v)I\big)$
\begin{equation}\label{pro2}
\pi_2(\v)=\frac{1}{2\pi
i}\oint_{|\zeta|=\delta}(\zeta-P_{\v}P_{\v}^\ast)^{-1}d\zeta,
\qquad \forall\v\in\tilde{U}.
\end{equation}
Since $V_1(\v_0)=\textrm{Ker}\ P_{\v_0}^\ast
P_{\v_0}=\textrm{Ker}\ P_{\v_0}$, we get that $\textrm{dim}\
V_1(\v_0)=1$, whence from Lemma 4.10 \cite{kato} we have that
$\textrm{dim}\ V_1(\v)=1$ for every $\v\in\tilde{U}$.\\
Therefore, if $\phi(\v_0; . )$ is an eigenfunction in $\textrm{Ker}\ P_{\v_0}=V_1(\v_0)$, then, provided one shrinks the neighborhood $\tilde{U}$, $\phi_1(\v; .
)=\pi_1(\v)\phi(\v_0; . )$ spans $V_1(\v)$ for any $\v\in\tilde{U}$. As a consequence of globally elliptic operators theory (see \cite{h3}) it turns out that $V_1(\v)\subset{\mathcal S}(\R)$; furthermore, since the operator $P_{\v}^\ast P_{\v}$ depends smoothly on $\v\in\tilde{U}$, the
same is true for the projector $\pi_1(\v)$, hence
$\phi_1(\v;.)\in C^\infty\big(\tilde{U},{\mathcal S}(\R)\big)$.\\
As for the second statement, by virtue of \eqref{autohom} it is
easily seen that
\begin{equation}\label{omII}
\phi_1(x',s\xi',s^{-\frac{1}{h+1}}t)=s^{\frac{1}{2(h+1)}} \phi_1(x',\xi',t),\,s>0.
\end{equation}
If $h$ is an odd integer, we easily check that $\phi_i(\varrho,-t)\in V_i(\varrho)$ ($i=1,2$); thus if we replace $\phi_i$ by the linear combination
$\phi_i(\varrho,t)\pm\phi_i(\varrho,-t)\in V_i(\varrho)$, that does not identically vanish  near $\varrho_0$, we get the parity in the $t-$variable we were seeking
for. Indeed the $\phi_1$ and $\phi_2$ have exactly the same parity because of the selfadjointness of the localized operator
$P_\varrho$ at $\varrho=\varrho_0$ (see Remark \ref{selfadj2} below); this yields $V_1(\varrho_0)=V_2(\varrho_0)$ whence one trivially gets $\phi_1(\varrho_0,t)=\beta\phi_2(\varrho_0,t)$ for some $0\neq\beta\in\C$.\\
Finally, in view of \eqref{omII} one obtains $\int
|\phi_i(x',s\xi',t)|^2 dt=\int |\phi_i(x',\xi',t)|^2 dt,$ so that
a standard normalization argument proves the last statement of
Lemma \ref{quad}.
\fine\\
If the localized operator is selfadjoint at $\v$, i.e. $P^\ast_\v=P_\v$, the eigenspace associated with the eigenvalue $\l_{j_0}$ is contained in
$V_1(\varrho)=V_2(\varrho)$ and both of these spaces are one-dimesional. As an immediate consequence, we have that $V_1(\varrho)=V_2(\varrho)=\textrm{Ker}\
(P_{\v}-\l_{j_0}(\v)I).$ Therefore, if $P$ is selfadjoint (and hence $P_\v$) the construction carried out in Lemma \ref{quad} is much easier.
\begin{remark}\label{selfadj}
Assume that $P^\ast_\v=P_\v$ for any $\v\in U$. Then one gets
that, for every $\v\in U$,
$$V_1(\varrho)=V_2(\varrho)=\textrm{Ker}\ (P_{\v}-\l_{j_0}(\v)I).$$
\end{remark}
Let us introduce a second remark which will be useful in Section \ref{mod}.
\begin{remark}\label{selfadj2}
Due to \eqref{avr}, \eqref{avreven}, \eqref{avrodd} one has that $\l_{j_0}(\v_0)=0$ and that $P^\ast_{\v_0}=P_{\v_0}$ for ${\v_0}=(x_0',\xi_0')$ whence
$$V_1(\varrho_0)=V_2(\varrho_0)=\textrm{Ker}\ (P_{\v_0}).$$
Therefore, we can choose an eigenfunction $\phi(\v;\cdot)$ in  $\textrm{Ker}\ (P_{\v_0})$ such that $\phi(\v_0;\cdot)=\phi_1(\v_0;\cdot)=\phi_2(\v_0;\cdot)$.
Arguing similarly as done in the proof of Lemma \ref{quad}, it is easily seen that $\phi(\v;\cdot)$ is smoothly dependent on the parameter $\v$.

\end{remark}

We can now take advantage of Lemma \ref{quad} to describe the map properties of $P_\v$.
\begin{lemma}\label{inv}
Denote by $V_j(\v)^\perp$ the $L^2(\R)$-orthogonal space of $V_j(\v)$ for $j=1,2$. For any $\v\in U$, the maps
$$
P_\v: B^2_h(\R)\cap V_1(\v)^\perp\longrightarrow V_2(\v)^\perp,\qquad P_\v: V_1(\v)\longrightarrow V_2(\v),
$$
are continuous isomorphisms. Moreover, if we define, for any $\v\in U$, the map
\begin{equation}\label{mappa}
\left\{\begin{array}{l}
E_\v: L^2(\R)\longrightarrow B^2_h(\R)\vspace{0,5cm}\\
E_\v:=\big(P_{\v\ |B^2_h(\R)\cap
V_1(\v)^\perp}\big)^{-1}\circ\big(I-\pi_2(\v)\big)
\end{array}\right.
\end{equation}
we have that
\begin{equation}\label{relaz}
P_\v E_\v+\pi_2(\v)=I,\quad E_\v P_\v+\pi_1(\v)=I,
\end{equation}
where, as above, $\pi_j(\v)$ represents the orthogonal projector
onto $V_j(\v)$ for $j=1,2$.
\end{lemma}
{\bf Proof.} A trivial check shows that
$$
P_\v^\ast(\zeta-P_\v P_\v^\ast)=(\zeta-P_\v^\ast
P_\v)P_\v^\ast,\quad\v\in\Z,
$$
thus, by applying the operators $(\zeta-P_\v^\ast P_\v)^{-1}$ and $(\zeta-P_\v P_\v^\ast)^{-1}$ to each side of the above equation, by \eqref{pro1} and \eqref{pro2}
we get $P_\v^\ast\pi_2(\v)=\pi_1(\v)P_\v^\ast$ for any $\v\in U$, whence $P_\v^\ast\big(V_2(\v)\big)\subset V_1(\v)$. As a consequence, for any $f\in B^2_h(\R)\cap
V_1(\v)^\perp$ and any $g\in V_2(\v)$ we get
$$\langle P_\v f, g\rangle_{L^2(\R)}= \langle f, P_\v^\ast g\rangle_{L^2(\R)}=0, $$
whence it follows that $P_{\v}(B^2_h(\R)\cap V_1(\v)^\perp)\subset V_2(\v)^\perp$. Similarly we see that $P_\v\big(V_1(\v)\big)\subset V_2(\v)$.\\
Since $\textrm{Ker}\, P_\v=\textrm{Ker}\, P_\v^\ast P_\v\subset V_1(\v)$, $P_{\v\ |B^2_h(\R)\cap V_1(\v)^\perp}$ is trivially injective; furthermore, by proceeding
as in Lemma 2.7 \cite{Parenti3}, the map turns out to be also surjective so that
\eqref{mappa} is well-defined and satisfies \eqref{relaz} by construction.\\
Finally, note that $\pi_2(\rho)$ is a smoothing global pseudodifferential operator in $\R$ (see Def. \ref{bbS} in Appendix or \cite{h3}), for its $(t,\tau)$-symbol
is given by
\begin{equation}\label{simbproi}
\sigma_{(t,\tau)}(\pi_2(\rho))=
e^{-it\cdot\tau}\phi_2(\rho;t)\hat{\overline{\phi_2}}(\rho; \tau).
\end{equation}
Since $P_\v$ is a $h-$globally elliptic pseudodifferential operator of order $2$ (see Def. \ref{bbS} and  Section 1.5 \cite{marco}), $P_\v$ has a $h-$globally
parametrix; therefore, from
the identity $P_\v E_\v=I-\pi_2(\v)$ it follows that $E_\rho$ is actually a $h-$globally pseudodifferential operator of order $-2$.\fine\\
 At last we are now in a
position to define the operators $h_{\phi_1}^+(\v), h_{\phi_2}^-(\v)$ and the corresponding localized system \eqref{locsis}.
\begin{definition}\label{h}
For any $\v\in U$, consider the following operators
$$
h_{\phi_1}^+(\v): L^2(\R)\longrightarrow \C,\quad
h_{\phi_1}^+(\v)f=\<f,\phi_1(\v,.)\>_{L^2(\R)},
$$
$$
h_{\phi_2}^-(\v): \C\longrightarrow V_2(\v),\qquad\quad
h_{\phi_2}^-(\v)\theta=\theta \phi_2(\v,.),
$$
and define the map
\begin{equation*}
{\cal A_\varrho}= \left(\begin{array}{cc}
\ds P_\varrho & h_{\phi_2}^-(\v)\\
h_{\phi_1}^+(\v) & 0
\end{array}\right):{\cal S}(\R)\times\C\longrightarrow {\cal S}(\R)\times\C
\end{equation*}
\end{definition}
As a consequence of the whole construction we get the following
theorem.
\begin{lemma}\label{inv2}
For every $\v\in U$, the localized system $A_\varrho$ is
invertible and its inverse is given by the map
\begin{equation*}
{\cal E_\varrho}= \left(\begin{array}{cc}
\ds E_\varrho & h_{\phi_1}^-(\v)\\
h_{\phi_2}^+(\v) & -\ell_{\frac{2}{h+1}}(\v)
\end{array}\right):{\cal S}(\R)\times\C\longrightarrow {\cal S}(\R)\times\C,
\end{equation*}
where $E_\v$ is defined in \eqref{mappa} and
$\ell_{\frac{2}{h+1}}(\v):=\<P_\v\phi_1(\v,.),\phi_2(\v,.)\>_{L^2(\R)}$
is a smooth function on $\v\in U$.
\end{lemma}
{\bf Proof.} By a direct computation the proof readily follows.\\
Defined the localized system ${\cal A}_\v, \;\v\in\Z$, we need to go back to the system
\eqref{sist} in all the variables. To this aim we have to define  the operators $H^\pm$, of which $h_{\phi_1}^+(\v), h_{\phi_2}^-(\v)$ are the corresponding
localized operators depending on the parameter $\v=(x',\xi')\in U$. Roughly speaking, $H^\pm$ are obtained by quantizing the functions $\phi_1(x',\xi',t),
\phi_2(x',\xi',t)$
 with respect to all the variables $(x',\xi';t,\tau)$.\\
In order to treat all these operators, together with their
localized operators, we need a pseudodifferential calculus
tailored to our anisotropic setting.
Since the calculus is very technical, we prefer to postpone this point to Appendix II, where the required classes of symbols and the corresponding composition rules are discussed.\\
Furthermore, in order to make the exposition more pleasant, we
carry out a local construction of the parametrix \eqref{param},
i.e. we assume that the neighborhood $U$ in Theorem \ref{inv} is
actually the whole cotangent space $T^\ast\R^{n-1}\setminus 0$.
The microlocal version is only a technical matter and we refer the
reader to Section 4 \cite{Parenti3} for the details.
\begin{definition}\label{hermite}
The Hermite operator $H^-$ is defined by
\begin{equation}\label{her}
\left\{\begin{array}{l}
H^-:C^\infty_0(\R_{x'}^{n-1})\longrightarrow C^\infty_0(\R_{t,x'}^{n})\vspace{0,5cm}\\
\ds (H^-f)(t,x')=(2\pi)^{-(n-1)}\int\int
e^{i(x',\xi')}\phi_2(x',\xi',t)\hat{f}(\xi')d\xi',
\end{array}\right.
\end{equation}
and the co-Hermite operator $H^+$ by
\begin{equation}\label{coher}
\left\{\begin{array}{l}
H^+:C^\infty_0(\R_{t,x'}^{n})\longrightarrow C^\infty(\R_{x'}^{n-1})\vspace{0,5cm}\\
\ds (H^+g)(x')=(2\pi)^{-(n-1)}\int\int e^{i(x',\xi')}\overline{\phi_1(x',\xi',t)}\hat{g}(t,\xi')d\xi'dt.
\end{array}\right.
\end{equation}
Accordingly with Definition \ref{def:H}  in Appendix II,
one has that $H^-\in \textrm{OP}H_h^\frac{1}{2(h+1)}$ and
$H^+\in\textrm{OP}H_h^{*\frac{1}{2(h+1)}}$. This is a consequence
of Lemma \ref{quad} and Remark \ref{ghherm} with
$m=\frac{1}{2(h+1)}$ and $j=0$.
\end{definition}
We can now state the main theorem of this section.
\begin{theorem}\label{main}
Consider the operator
\begin{equation}
{\cal A}= \left(\begin{array}{cc}
\ds P & H^-\\
H^+ & 0
\end{array}\right).
\end{equation}
There exist
\begin{enumerate}
\item $\ds E\in \textrm{OP}S^{-\frac{2}{h+1}}_{\frac{1}{h+1},\frac{1}{h+1}}$ $($more precisely, $E\in \textrm{OP}S^{-2,-2}_h$ with asymptotic expansion
$\sigma(E)\sim\sum_{j\geq 0} {e}_{-j}(t,x',\tau,\xi'),\; e_{-j}\in S_h^{-2,-2+j/h}$, see Definition \ref{definition11} and the notes below in Appendix
\ref{calc}$)$;
\item $K^-\in \textrm{OP}H_h^\frac{1}{2(h+1)}$ with symbol $\sigma(K^-)\sim\sum_{j\geq 0} \psi'_{-j/(h+1)}(x',\xi',t)$, where
$\psi'_{0}(x',\xi',t):=\phi_1(x',\xi',t);$
\item $K^+\in \textrm{OP}H_h^{*\frac{1}{2(h+1)}}$ with symbol $\sigma(K^+)\sim\sum_{j\geq 0} \psi_{-j/(h+1)}(x',\xi',t)$,
where $\psi_{0}(x',\xi',t):=\phi_2(x',\xi',t);$
\item a pseudodifferential operator $L=\ell(x',D_{x'})\in\textrm{OP}S_{1,0}^{\frac{2}{h+1}}(\R^{n-1})$, with
(poly)homogeneous asymptotic expansion $\ds \ell(x',\xi')\sim\sum_{j\geq 0} \ell_{\frac{2}{h+1}-\frac{j}{h+1}}(x',\xi')$ (see Def. 18.1.5 Vol.III \cite{ho}),
and principal symbol given by $\ell_{\frac{2}{h+1}}(x',\xi')=\<P_{(x',\xi')}\phi_1(x',\xi',.),\phi_2(x',\xi',.)\>_{L^2(\R)}$,
\end{enumerate}
such that the system
\begin{equation}\label{param2}
{\cal E}= \left(\begin{array}{cc}
\ds E & K^-\\
K^+ & -L
\end{array}\right)
\end{equation}
is a two-sided parametrix of $\cal A$; namely, ${\cal AE}-I$ and
${\cal EA}-I$ map
${\cal E}'(\R^n)\times{\cal E}'(\R^{n-1};\C)\longrightarrow C^\infty(\R^n)\times C^\infty(\R^{n-1};\C).$\\
Furthermore, $P$ is hypoelliptic with loss of $2h/(h+1)+\delta$ derivatives if and only if so is $L$ with loss $\delta$ $($see \eqref{defipo}$)$.
\end{theorem}
{\bf Proof.} The construction of the parametrix $\cal E$ requires
the full strength of the calculus developed in Appendix II
and is very involved. The strategy is to look for a system of the
form \eqref{param2} such that
\begin{equation}\label{primo}
{\cal AE}-I\equiv 0,
\end{equation}
where ${\cal AE}-I$ is a smoothing operator. A similar argument applied to $\cal A^*$ yields ${\cal A}^*\tilde{\cal E}-I\equiv 0$, whence one gets $\tilde{\cal E}^*{\cal A}-I\equiv 0$; by standard calculus it turns out that ${\cal E - \tilde{E}^*}\equiv 0$ so that $\cal E$ is also a left parametrix of $\cal A$.\\
The equation \eqref{primo} amounts to saying that
\begin{equation}\label{equazioni}
 \underbrace{PE+H^-K^+-I\equiv 0}_{(i)},\quad \underbrace{H^+E\equiv 0}_{(ii)},\quad \underbrace{PK^--H^-L\equiv 0}_{(iii)},\quad \underbrace{H^+K^--I\equiv 0}_{(iv)}.
\end{equation}
Therefore, in order to get the parametrix $\cal E$ we have to look for operators $E, K^\pm, L$ solving these equations.\\
In particular, it is enough to solve the equation \eqref{equazioni} $(i)$ modulo $\textrm{OP}S_h^{0,\infty}$ (see Proposition 3.1.2 \cite{h2}). To take advantage of
that observation we need to introduce the localized operators $P_\v^{(2+r)}$ of any order, i.e. for  $r=0,1,2,...,$
\begin{equation}\label{lochigh}
\ds P_\v^{(2+r)}= {\sum_{{\scriptsize\begin{array}{c}
        \a,\b \in \mathbb{Z}_+ \\
        \ds \frac{\a}{h}+\b+\frac{h+1}{h}j=2+\frac{r}{h}
        \end{array}}}}
 \frac{1}{\a ! \b !}(\partial^\a_{x_1}\partial^\b_{\xi_1}p_{2-j})(\varrho)t^\a D_t^\b,
\end{equation}
where $p_{2-j}$ denotes the homogeneous term of order $2-j$ in the asymptotic expansion of the complete symbol of $P$ in \eqref{grushin}.
A straightforward computation shows that $P_\v^{(2)}$ coincides exactly with the usual localized operator $P_\v$ in \eqref{locgrushin}. \\
We point out that $P_\v^{(2+r)}$ is the classical quantization of the Taylor expansion of the symbol of $P$ at the order $2+r/h$ near $\Z$; hence, the symbol of
$P-\sum_{r\geq 0}P_\v^{(2+r)}$ vanishes of infinite order on $\Sigma$ so that $P-\sum_{r\geq 0}P_\v^{(2+r)}\in\textrm{OP}S_h^{2,\infty}$. Therefore, the equation
\eqref{equazioni} $(i)$ can be replaced by
\begin{equation}\label{equazioni2}
 \big(\sum_{r\geq 0}P_\v^{(2+r)}\big)E+H^-K^+-I=0\quad \textrm{mod.}\,\textrm{OP}S_h^{0,\infty}\qquad\qquad\qquad\qquad\qquad\qquad(i')
\end{equation}
Finally, denoting by $\sigma(E)\sim\sum_{j\geq 0}
e_{-j}(t,x',\tau,\xi'),\; e_{-j}\in S_h^{-2,-2+j/h}$ the
asymptotic expansion of the operator $E$ (we are seeking for), we
define
$$
E_{\v=(x',\xi')}^{(-2-j)}:=e_{-j}(t,x',D_t,\xi'),\;j\geq 0,
$$
as the pseudodifferential operator obtained by quantizing
$e_{-j}(t,x',\tau,\xi')$ only in the variables $(t,\tau)$.
According to this notation, we choose
$E_{\v=(x',\xi')}^{(-2)}:=E_\v$ where $E_\v$
is the pseudodifferential operator defined in Lemma \ref{inv}. \\
By applying the composition rules in Appendix II, we can
rephrase the equations \eqref{equazioni2} $(i')$ and
\eqref{equazioni} $(ii)-(iv)$ in the following ``algebraic"
relations, which have to be solved  at each degree  $s=0,1,...$ of
homogeneity
\begin{itemize}
\item[(I)]
$\ds\sum_{(h+1)|\a|+j+q=s}\frac{1}{\alpha!i^{|\a|}}\sigma\big(
\partial^\a_{\xi'}P_{(x',\xi')}^{(2+j)}\#
\partial^\a_{x'}E^{(-2-q)}_{(x',\xi')}\big)$
$$+\sum_{(h+1)|\a|+q=s}\frac{1}{\alpha!i^{|\a|}} e^{-it\cdot\tau}\partial_{\xi'}^\alpha\phi_{2}(x',\xi';t) \partial_{x'}^\a\widehat{\overline{\psi}}_{-q/(h+1)}(x',\xi';\tau)=\left\{\begin{array}{l}
                                                         1\;\;\tr{if}\;s=0\\
                                                         0\;\;\tr{if}\;s\geq1
                                                        \end{array}\right.\\
 $$
\item[(II)]$\ds \sum_{(h+1)|\a|+j=s}\frac{1}{\alpha!i^{|\a|}}
\overline{\Big(\partial^\a_{\xi'}E^{(-2-j)}_{(x',\xi')}\Big)^*
\big(\partial_{x'}^\alpha\phi_{1}(x',\xi';\cdot)\big)}=0$
\item[(III)]$\ds\sum_{(h+1)|\a|+j+q=s}\frac{1}{\alpha!i^{|\a|}}
\partial^\a_{\xi'}P_{(x',\xi')}^{(2+j)} \big(
\partial^\a_{x'}\psi'_{-q/(h+1)}(x',\xi';\cdot)\big)$
$$-\sum_{(h+1)|\a|+q=s}\frac{1}{\alpha!i^{|\a|}} \partial_{\xi'}^\alpha\phi_{2}(x',\xi';t) \partial_{x'}^\a \ell_{\frac{2}{h+1}-\frac{q}{h+1}}(x',\xi')=0;\\
 $$
\item[(IV)] $\ds\sum_{(h+1)|\alpha|+q=s}\frac{1}{\alpha
!i^{|\alpha|}}
\Big(\partial_{x'}^\alpha{\psi'}_{-q/(h+1)}(x',\xi';\cdot),\partial_{\xi'}^\alpha{\phi}_{1}
(x',\xi';\cdot)\Bigr)_{L^2(\R)} =\left\{\begin{array}{lll}
1,\,\,\mathrm{if}\,\, s=0\\
0,\,\,\mathrm{if}\,\,s\geq 1.\end{array}\right.$
\end{itemize}
We have thus to determine all the symbols $e_{s}, \psi_{-s/(h+1)},\,\,{\psi'}_{-s/(h+1)}, \ell_{2/(h+1)-s/(h+1)}.$\\
This crucial task will be accomplished by using an inductive
procedure. We are now going to show that the first step of the
iteration (i.e. $s=0$) is an immediate consequence of the
construction carried out in this section. Precisely, as for $s=0$,
we have already chosen $E_{\v=(x',\xi')}^{(-2)}:=E_\v,\,
\psi_{0}(x',\xi',t):=\phi_2(x',\xi',t),\,
\psi'_{0}(x',\xi',t):=\phi_1(x',\xi',t),\,
\ell_{\frac{2}{h+1}}(x',\xi'):=\<P_{(x',\xi')}\phi_1(x',\xi',.),\phi_2(x',\xi',.)\>_{L^2(\R)}$
and
 we are thus left to check that
\begin{itemize}
\item[$\tr{(I)}_{s=0}$] $\ds\quad\sigma\big( P_{(x',\xi')}^{(2)}\#
E^{(-2)}_{(x',\xi')}\big)+ e^{-it\cdot\tau}\phi_{2}(x',\xi';t) \widehat{\phi}_{2}(x',\xi';\tau)=1$ or equivalently, by \eqref{simbproi}\\
$$\ds\hspace{-8cm}\sigma\big( P_{(x',\xi')}^{(2)}\#
E^{(-2)}_{(x',\xi')}+\pi_2(x',\xi')\big)=1$$
\item[$\tr{(II)}_{s=0}$]$\ds\quad
\overline{\Big(E^{(-2)}_{(x',\xi')}\Big)^*
\big(\phi_{1}(x',\xi';\cdot)\big)}=0$
\item[$\tr{(III)}_{s=0}$]$\ds\quad P_{(x',\xi')}^{(2)}
\big(\phi_{1}(x',\xi';\cdot)\big)(t)=
\<P_{(x',\xi')}\phi_1(x',\xi',.),\phi_2(x',\xi',.)\>_{L^2(\R)}\phi_{2}(x',\xi';t)
;$ \item[$\tr{(IV)}_{s=0}$] $\ds\quad
\Big({\phi}_{1}(x',\xi';\cdot),{\phi}_{1}
(x',\xi';\cdot)\Bigr)_{L^2(\R)} =1$
\end{itemize}
The above relations follow immediately from Lemma \ref{inv2},
since ${\cal E}_\v$ is the (right) inverse of the matrix operator
${\cal A}_\v$. At this point the iterative procedure can start
and, arguing similarly as done in Theorem 4.1 \cite{Parenti3} or also in $(3-6)$ \cite{bmt2}, we
can solve the equations $\tr{(I)}-\tr{(IV)}$ in $e_{-s},
\psi_{-s/(h+1)},\,\,{\psi'}_{-s/(h+1)},$
$\ell_{2/(h+1)-s/(h+1)}$ for any $s\geq1$.\\
In order to state our results about hypoellipticity it is
essential to know explicitly the operator $L=\ell(x',D_{x'})$. For
this reason we complete the first part of the proof
by writing out the symbols $\ell_{2/(h+1)-s/(h+1)},\;s\geq1$ of its asymptotic expansion, as they come out from the iteration.\\
Precisely, assume that ${\psi'}_{-q/(h+1)}$ and
$\ell_{2/(h+1)-q/(h+1)}$, $0\leq q<s$, have already been
determined, we perform the following decomposition
$${\psi'}_{-s/(h+1)}(x',\xi';\cdot)={\psi'}_{-s/(h+1),1}(x',\xi';\cdot)+
{\psi'}_{-s/(h+1),2}(x',\xi';\cdot)\in V_1(x',\xi')\oplus
V_1(x',\xi')^\perp,$$ Since ${\psi'}_{-s/(h+1),1}=
\big({\psi'}_{-s/(h+1)},{\phi}_1\big)_{L^2(\R)}{\phi}_1$,\;
${\psi'}_{-s/(h+1),1}$ is uniquely determined by (IV) above, so
that one gets
\begin{equation}\label{psi'1}
\Bigl({\psi'}_{-s/(h+1),1}(x',\xi';\cdot),{\phi}_{1}(x',\xi';\cdot)\Bigr)
_{L^2(\R)}=\hspace{3cm}
\end{equation}
$$=-\sum_{\text{\scriptsize$\begin{matrix}(h+1)|\alpha|+q=s\\ 0\leq q<s\end{matrix}$}}
\frac{1}{\alpha
!i^{|\alpha|}}\Bigl(\p_{x'}^\alpha{\psi'}_{-q/(h+1)}(x',\xi';\cdot),
\p_{\xi'}^\alpha{\phi}_{1}(x',\xi';\cdot)\Bigr)_{L^2(\R)}.$$ As
for $\ell_{2/(h+1)-s/(h+1)},$  by taking in (III) above the
$L^2(\R)$-scalar product with $\phi_2$ and by Lemma \ref{inv} one
has
\begin{equation}
\ell_{2/(h+1)-s/(h+1)}(x',\xi')=\Bigl(P_{(x',\xi')}^{(2)}{\psi'}_{-s/(h+1),1}
(x',\xi';\cdot),\phi_2(x',\xi';\cdot)\Bigr)_{L^2(\R)}+\hspace{1cm}
\label{ell}\end{equation}
$$-\sum_{\text{\scriptsize$\begin{matrix}(h+1)|\alpha|+q=s\\ 0\leq q<s\end{matrix}$}}\hspace{-.2cm}
\frac{1}{\alpha
!i^{|\alpha|}}\p_{x'}^\alpha\ell_{2/(h+1)-q/(h+1)}(x',\xi')\Bigl(
\p_{\xi'}^\alpha\phi_{2}(x',\xi';\cdot),\phi_2(x',\xi';\cdot)\Bigr)_{L^2(\R)}\hspace{-.3cm}+$$
$$+\sum_{\text{\scriptsize$\begin{matrix}(h+1)|\alpha|+j+q=s\\ 0\leq q<s\end{matrix}$}}
\frac{1}{\alpha !i^{|\alpha|}}\Bigl(\p_{\xi'}^\alpha
P_{(x',\xi')}^{(2+j)}\p_{x'}^\alpha{\psi'}_{-q/(h+1)}
(x',\xi';\cdot),\phi_2(x',\xi';\cdot)\Bigr)_{L^2(\R)}.$$ Finally,
as for ${\psi'}_{-s/(h+1),2},$ again from (III) above we have
\begin{equation}\label{psi'2}
{\psi'}_{-s/(h+1),2}(x',\xi';\cdot)=E^{(-2)}_{(x',\xi')}\Biggl[
-P_{(x',\xi')}^{(2)}{\psi'}_{-s/(h+1),1}(x',\xi';\cdot)+
\end{equation}
$$
\sum_{\text{\scriptsize$\begin{matrix}(h+1)|\alpha|+q=s\end{matrix}$}}
\frac{1}{\alpha
!i^{|\alpha|}}\p_{\xi'}^\alpha\phi_{2}(x',\xi';\cdot)\p_{x'}^\alpha
\ell_{2/(h+1)-q/(h+1)}(x',\xi')+$$
$$-\sum_{\text{\scriptsize$\begin{matrix}(h+1)|\alpha|+j+q=s\\ 0\leq q<s\end{matrix}$}}
\frac{1}{\alpha !i^{|\alpha|}}\p_{\xi'}^\alpha
P_{(x',\xi')}^{(2+j)}\p_{x'}^\alpha{\psi'}_{-q/(h+1)}(x',\xi';\cdot)
\Biggr].$$
Note that the argument of $E^{(-2)}_{(x',\xi')}$ in the right hand side of the identity \eqref{psi'2} belongs to $V_2((x',\xi'))^\perp$, because of \eqref{ell}.\\
Let us complete the proof by showing that $P$ is hypoelliptic with loss of $2h/(h+1)+\delta$ derivatives if and only if so is $L$ with loss $\delta$.\\
Taking into account \eqref{defipo}, suppose that $P$ satisfies $(\textrm{He})_{2h/(h+1)+\delta}$ (see \eqref{defipo}), we are going to show that $L$ verifies $(\textrm{He})_{\delta}$.\\
For the sake of simplicity, we assume that $\Omega=\R^{n-1}$.
Given any $f\in{\cal D}'(\R^{n-1})$ we have thus to prove that if
$Lf\in H^s_{loc}(\R^{n-1})$ then $f\in
H^{s+2/(h+1)-\delta}_{loc}(\R^{n-1})$. By Lemma \ref{Sobolev} we
get $H^-Lf\in H^{s}_{loc}(\R^{n})$ and by $(iii)$ of
\eqref{equazioni} we have $PK^-f-H^-Lf\in C^\infty(\R^{n-1})$,
whence $PK^-f\in H^{s}_{loc}(\R^{n})$. In view of the
hypoellipticity of $P$ we obtain $K^-f\in
H^{s+2/(h+1)-\delta}_{loc}(\R^{n})$ and, again by Lemma
\ref{Sobolev},
$H^+K^-f\in H^{s+2/(h+1)-\delta}_{loc}(\R^{n-1})$. Finally, from $(iv)$ of \eqref{equazioni} we conclude that $f\in H^{s+2/(h+1)-\delta}_{loc}(\R^{n-1})$.\\
On the other hand, given any $g\in{\cal D}'(\R^{n})$, assume now
that $L$ verifies $(\textrm{He})_{\delta}$, and $Pg\in
H^{s}_{loc}(\R^{n})$. Since $\cal E$ is also a left parametrix of
$\cal A$, i.e. ${\cal AE}-I\equiv 0$, we have:
\begin{equation}\label{sinistra}
EP+K^-H^+=I\qquad K^+P-LH^+=0.
\end{equation}
Lemma \ref{Sobolev} yields $EPg\in H^{s+2/(h+1)}_{loc}(\R^{n})$ and the first equation in \eqref{sinistra} gives $EPg+K^-H^+g-g\in C^\infty(\R^n)$, whence $K^-H^+g-g\in H^{s+2/(h+1)}_{loc}(\R^{n})$. Moreover, by Lemma \ref{Sobolev} one has $K^+Pg\in H^s_{loc}(\R^{n-1})$ and from the second equation in \eqref{sinistra} it follows that $K^+Pg-L H^+g\in C^\infty(\R^{n-1})$, therefore $L H^+g\in H^s_{loc}(\R^{n-1})$. Since $L$ is $(\textrm{He})_{\delta}$, one has $H^+g\in H^{s+2/(h+1)-\delta}_{loc}(\R^{n-1})$
and, again by Lemma \ref{Sobolev}, $K^-H^+g\in H^{s+2/(h+1)-\delta}_{loc}(\R^{n})$, whence we finally get $g\in H^{s+2/(h+1)-\delta}_{loc}(\R^{n})$.\fine\\
In view of Theorem \ref{main} we are thus reduced to studying the h.e. with minimal loss of derivatives of the operator $L=\ell(x',D_{x'})$ in $(n-1)-$variables,
with non classical asymptotic expansion $\ell(x',\xi')\sim\sum_{j\geq 0} \ell_{\frac{2}{h+1}-\frac{j}{h+1}}(x',\xi')$.\\
The analysis of such $h$-homogeneous operators is carried out in
Appendix I, where the h.e. of  ``general" operators
${\rm Op}(a)$, $a\sim\sum_{j\geq 0} a_{m'-\frac{j}{h+1}}$ is
studied. The results in Appendix I here apply with
$\nu=n-1$ and $m'=2/(h+1)$, and Propositions \ref{cns} and
\ref{h=3} show that only the first three terms
$\ell_{\frac{2}{h+1}-\frac{j}{h+1}},\;j=0,1,2$ in the asymptotic
expansion of $\ell(x',D_{x'})$ really matter in the analysis of the h.e. with minimal loss of regularity.
For this reason, in
the following remarks  we compute explicitly the dependence of
such terms from the operator $P$.
\begin{remark}\label{struttura}
From \eqref{ell} we get
\begin{eqnarray}
\ell_{\frac{2}{h+1}}(x',\xi')&=&
\<P^{(2)}_{(x',\xi')}\phi_1,\phi_2\>_{L^2(\R_t)},\label{simbc1}\\
\ell_{\frac{2}{h+1}-\frac{1}{h+1}}(x',\xi')&=&
\<P^{(3)}_{(x',\xi')}\phi_1,\phi_2\>_{L^2(\R_t)},\label{simbc2}\\
\ell_{\frac{2}{h+1}-\frac{2}{h+1}}(x',\xi')&=&
\<P^{(4)}_{(x',\xi')}\phi_1,\phi_2\>_{L^2(\R_t)}\label{simbc3}\\
& & -\<E^{(-2)}_{(x',\xi')}\big(P^{(3)}_{(x',\xi')}\phi_1\big),
\big(P^{(3)}_{(x',\xi')}\big)^\ast\phi_2\>_{L^2(\R_t)},\nonumber
\end{eqnarray}
 where $h>1$  in the last equation.\\
If $h$ is an odd integer, the operator $P^{(3)}_{(x',\xi')}$ (see
\eqref{lochigh}) flips the parity in the $t-$variable; therefore,
from 2. of Lemma \ref{quad}, we immediately get that
$\ell_{\frac{2}{h+1}-\frac{1}{h+1}}(x',\xi')\equiv 0.$\\
This is actually the reason for which the loss of regularity of
$P$ strictly depends on the parity of $h$.
\end{remark}

Furthermore, if the localized operator $P_\v=P^{(2)}_{(x',\xi')}$
in \eqref{locgrushin} is selfadjoint (i.e., $b_1(0,x',\xi')$ is a
real-valued function), then by Remark \ref{selfadj} the functions
$\phi_1,\phi_2$ can be chosen equal to the normalized eigenfuction
$\phi(t,x',\xi')$ associated with the eigenvalue
$\l_{j_0}(x',\xi')$ of $P_\v$. Therefore, we have the following.
\begin{remark}\label{struttura2}
Again from \eqref{ell} we get (here $h>1$)
\begin{eqnarray}
\ell_{\frac{2}{h+1}}(x',\xi')&=&
\<P^{(2)}_{(x',\xi')}\phi,\phi\>_{L^2(\R_t)}=\l_{j_0}(x',\xi'),\label{simbc1}\\
\ell_{\frac{2}{h+1}-\frac{1}{h+1}}(x',\xi')&=&
\<P^{(3)}_{(x',\xi')}\phi,\phi\>_{L^2(\R_t)},\label{simbc2}\\
\ell_{\frac{2}{h+1}-\frac{2}{h+1}}(x',\xi')&=&
\<P^{(4)}_{(x',\xi')}\phi,\phi\>_{L^2(\R_t)}\label{simbc3}\\
& &-\<E^{(-2)}_{(x',\xi')}\big(P^{(3)}_{(x',\xi')}\phi\big),
\big(P^{(3)}_{(x',\xi')}\big)^\ast\phi\>_{L^2(\R_t)},\nonumber
\end{eqnarray}
 with $h>1$  in the last equation.\\
Again, if $h$ is an odd integer, we have that
$\ell_{\frac{2}{h+1}-\frac{1}{h+1}}(x',\xi')\equiv 0.$
\end{remark}
We are now in a position to state the main theorem of this
Section.

\begin{theorem}\label{teoprinc}
Let $P$ be the operator \eqref{grushin} and assume that
\eqref{avr} holds (i.e., \eqref{avreven} or
\eqref{avrodd} according to the parity of $h$).  Let
$L=\ell(x',D_{x'})$, with $\ell\sim\sum_{j\geq 0}
\ell_{\frac{2}{h+1}-\frac{j}{h+1}}$, be the operator defined in
$4.$ of Theorem \ref{main}. Then $P$ can be
$C^\infty$-hypoelliptic only with loss of
$\delta\geq\frac{2h}{h+1}+\frac{1}{h+1}$ derivatives if $h$ is
even and of $\delta\geq2$ derivatives if $h$ is odd. Furthermore,
$\delta$ attains the above lower bounds (i.e., $P$ is h.e. with minimal
loss of derivatives) if and only if
\begin{itemize}
\item[$(I)$] $($$h$ even$)$:
$\ell_{\frac{1}{h+1}}(x'_0,\xi'_0)\neq 0$, and there exist a conic
neighborhood $U$ of $(x'_0,\xi'_0)$ in $T^\ast\mathbb{S}^{n-1}$
and a constant $c=c(U)>0$ such that, for every $(x',\xi')\in U$,
$$-{\rm{Re}}\  \ell_{\frac{2}{h+1}}(x',\xi')\ell_{\frac{1}{h+1}}(x',\xi')\leq
|\ell_{\frac{2}{h+1}}(x',\xi')|\sqrt{|\ell_{\frac{1}{h+1}}(x',\xi')|^2-c}\;;$$

\item[$(II)$] $($$h$ odd, $h\neq 3$$)$: $\ell_{0}(x'_0,\xi'_0)\neq
0$, and there exist a conic neighborhood $U$ of $(x'_0,\xi'_0)$ in
$T^\ast\mathbb{S}^{n-1}$ and a constant $c=c(U)>0$ such that, for
every $(x',\xi')\in U$,
$$-{\rm{Re}}\  \ell_{\frac{2}{h+1}}(x',\xi')\ell_{0}(x',\xi')\leq
|\ell_{\frac{2}{h+1}}(x',\xi')|\sqrt{|\ell_{0}(x',\xi')|^2-c}\;.$$
\end{itemize}
 Assume that $h=3$. Then $P$ satisfies the minimal loss of derivatives if\linebreak   $|\ell_{0}(x'_0,\xi'_0)|^2+\frac{1}{i}\{\overline{\ell_{\frac{1}{2}}},
    \ell_{\frac{1}{2}}\}(x'_0,\xi'_0)>0$, and there exist a conic neighborhood $U$ of $(x'_0,\xi'_0)$ in $T^\ast\mathbb{S}^{n-1}$ and a constant $c=c(U)>0$ such that, for every $(x',\xi')\in U$,
$$-{\rm{Re}}\  \ell_{\frac{1}{2}}(x',\xi')\ell_{0}(x',\xi')\leq
|\ell_{\frac{1}{2}}(x',\xi')|\sqrt{|\ell_{0}(x',\xi')|^2+\frac{1}{i}\{ \overline{\ell_{\frac{1}{2}}},\ell_{\frac{1}{2}}\}(x',\xi')-c}\;.$$ Finally, if the
$\ell_{\frac{1}{2}}(x',\xi')$ vanishes identically on a conic neighborhood of $(x'_0,\xi'_0)$, then the  condition $\ell_{0}(x'_0,\xi'_0)\neq 0$ is actually
necessary and sufficient.
\end{theorem}
The proof is a straightforward consequence of Theorem \ref{main} and Lemma \ref{loss} , Propositions \ref{cns} and \ref{h=3} in Appendix I.
Precisely, from the former the analysis of the hypoellipticity of
$P$ is reduced to the study of the hypoellipticity of
$L=\ell(x',D_{x'})$, which is carried out in the above mentioned
results in the appendix.

%Alla fine di tutto si fanno due casi: 1) in cui nei punti degeneri il loc è autoaggiunto e si usa per grushin 2) Caso operatori autoaggiunti che si usa per Op di kohn
\section{Kohn and Gilioli-Treves operators}\label{ko}
In this section we discuss the hypoelliticity of several examples. The main difficulty in doing that relies the lack of a complete description of the spectrum and
of the eigenfunctions of the anharmonic oscillator \eqref{locgrushin}. However we are going to show that in many cases this knowledge is not actually strictly
necessary for our aims.
\subsection{A Gilioli-Treves model}
Consider the following operator in $\R^2$:
\begin{equation*}
P=D_{x_1}^2+ax_1^{2h}D_{x_2}^2+\b(x_1)x_1^{h-1}D_{x_2}
\end{equation*}
where $0<a\in\R$, $h\in\N$ and $\b(x_1)\in C^\infty(\R)$ is a real
function. If its localized operator \eqref{locgrushin}
\begin{equation}\label{locex1}
P_\varrho=D_t^2 +a|\xi_2|^2t^{2h}+\b(0)\xi_2 t^{h-1}
\end{equation}
is injective at every point $\varrho=(x_2,\xi_2)$, then $P$ is
classically hypoelliptic with loss of $2h/(h+1)$. Suppose now that
$P_\varrho$ is not injective at some $\varrho_0=(x^0_2,\xi^0_2)$
and, firstly, assume that $h$ is an even integer. This means,
accordingly to \eqref{avreven}, that
$$|\b(0)|=\sqrt{a}(h+1)(2j_0+1),$$
for some $j_0\in\N$. In view of \eqref{avr} this implies that the
eigenvalue $\l_j(\varrho)$ of $P_\varrho$ vanishes at
$\varrho=\varrho_0$. Since $P$ is actually independent of the $x_2$-variable, a standard scaling
argument shows that $P_\varrho$ is never injective at any
$\varrho=(x_2,\xi_2)\in \textrm{Char} P$, so that
$$\l_{j_0}(x_2,\xi_2)\equiv 0,\qquad\forall(x_2,\xi_2)\in\R^2,\; \xi_2\neq0.$$
Moreover, since $P_\varrho^*=P_\varrho$, from Remark \ref{struttura2} and \eqref{lochigh} with $r=1$ we
have that
\begin{eqnarray*}
\ell_{\frac{2}{h+1}}(x_2,\xi_2)&=&\l_{j_0}(x_2,\xi_2)\equiv 0,\\
\ell_{\frac{1}{h+1}}(x_2,\xi_2)&=&
\<P^{(3)}_{(x_2,\xi_2)}\phi,\phi\>_{L^2(\R_{x_1})}
=\b'(0)\xi_2\<x_1^h\phi,\phi\>_{L^2(\R_{x_1})}=\b'(0)\xi_2
\|x_1^{h/2}\phi\|^2_{L^2(\R_{x_1})},
\end{eqnarray*}
where $\phi=\phi(x_1,x_2,\xi_2)$ is a normalized eigenfuction associated with the eigenvalue $\l_{j_0}(x_2,\xi_2)$ of $P_\v$. As a consequence of $(I)$ in Theorem \ref{teoprinc}, $P$ is $C^\infty-$hypoelliptic with minimal loss of $\frac{2h}{h+1}+\frac{1}{h+1}$ derivates if and only if $\b'(0)\neq 0$.\\
If $h$ is an odd integer, the analysis is more delicate. From
\eqref{avrodd} we have that
\begin{equation}\label{esempioodd}
\b(0)\xi_2=(-1)^{j_0}\sqrt{a}|\xi^0_2|-\sqrt{a}(h+1)(j_0+\theta(j_0))|\xi^0_2|
\end{equation}
and arguing as before yields
\begin{eqnarray*}
\ell_{\frac{2}{h+1}}(x_2,\xi_2)&=&
\l_{j_0}(x_2,\xi_2)\equiv 0,\\
\ell_{\frac{1}{h+1}}(x_2,\xi_2)&=&
\<P^{(3)}_{(x_2,\xi_2)}\phi,\phi\>_{L^2(\R_t)}\equiv 0\\
\ell_{0}(x_2,\xi_2)&=&
\<P^{(4)}_{(x_2,\xi_2)}\phi,\phi\>_{L^2(\R_t)}
 -\<E^{(-2)}_{(x_2,\xi_2)}\big(P^{(3)}_{(x_2,\xi_2)}\phi\big),
\big(P^{(3)}_{(x_2,\xi_2)}\big)^\ast\phi\>_{L^2(\R_t)}\\
&=& \frac{1}{2}\b''(0)\xi_2\|x_1^{(h+1)/2}\phi\|^2_{L^2(\R_{x_1})}
-(\b'(0)\xi_2)^2\<E^{(-2)}_{(x_2,\xi_2)}\big(x_1^h\phi\big),
\big(x_1^h\phi\big)\>_{L^2(\R_t)}.
\end{eqnarray*}
Due to Theorem \ref{teoprinc}, $P$ turns out to be hypoelliptic
with loss of $2$ derivatives if and only if
$\ell_{0}(x_2,\xi_2)\neq 0$; therefore, a trivial necessary
condition is that
\begin{equation}\label{esempiocond}
(\b'(0),\b''(0))\neq (0,0).
\end{equation}
Unlike the even case, this condition is not generally sufficient;
to this aim let us assume that $j_0=0$ in \eqref{esempioodd} so
that $\b(0)\xi_2=-\sqrt{a}h|\xi^0_2|$. According to \eqref{avr}
and \eqref{ord}, the first eigenvalue is identically zero so that
the localized operator $P^{(2)}_{(x_2,\xi_2)}$ is a
non negative operator as well its ``partial inverse''
$E^{(-2)}_{(x_2,\xi_2)}$. Therefore, if
$$\left\{\begin{array}{l}
\b'(0)= 0,\\
\b''(0)\neq 0,
\end{array}\right.\qquad \textrm{or}\qquad
                                      \b'(0)\neq 0 \Longrightarrow \left\{\begin{array}{l}
                                      \b''(0)\leq 0\quad \textrm{if}\quad \xi_2^0>0,\\
                                      \b''(0)\geq 0 \quad\textrm{if}\quad \xi_2^0<0,
                                      \end{array}\right.$$
then  $\ell_{0}(x_2,\xi_2)\neq 0$, whence the hypoellipticity of $P$ immediately follows.
\subsection{An extension of the Kohn operator:}\label{koext}
Let us now consider in $\R^2$ the following vector field
$$L=D_{x_1}+i g(x_1)D_{x_2}$$
where $g(x_1)$ is a real polynomial (or a real analytic function). Note that $L$ can be regarded as a generalization of the Mizohata operator $M=D_{x_1}-i x_1^h
D_{x_2}$ ($h\in\N$). Consider the sum of squares of complex vector fields
\begin{equation}\label{kext}
P=LL^*+\big(f(x_1)L\big)^*f(x_1)L
\end{equation}
$$=(1+f(x_1)^2)\big(D_{x_1}^2+g(x_1)^2D_{x_2}^2\big)-(1-f(x_1)^2)g'(x_1)D_{x_2}-2if'(x_1)f(x_1)L,$$
%$$P=LL^*=D_{x_1}^2+g(x_1)^2D_{x_2}^2-g'(x_1)D_{x_2}.$$
where $f(x_1)$ denotes any real analytic function. Its
characteristic set is given by $\Sigma=\{(x,\xi) | \xi_1=0,
g(x_1)=0, \xi_2\neq 0\}$. Therefore the hypoellipticity of $P$
strictly depends on the behaviour of $g$ near its real roots. Let
us begin our analysis in a small vertical strip centered in one of
these roots $x_1^0$, say $S_0=\{x\in\R^2\, |\, |x_1-x_1^0|<\e\}$. Without loss of generality, we can assume
that $x_1^0=0$ and we have that $g(x_1)=a(x_1)x_1^h$, where $h$
denotes the order of the root $x_1^0=0$ and $a(x_1)\in
C^\omega(\R), a(0)>0$ (similarly, if $a(0)<0$). As a consequence, one has that $\Sigma\cap T^\ast S_0=\{\xi_1=0=x_1, \, \xi_2\neq 0\}$ and
the operator $Q$ can be written as
$$P=(1+f(x_1)^2)\big(D_1^2+a(x_1)^2x_1^{2h}D_2^2\big)$$
$$-(1-f(x_1)^2)\big(ha(x_1)x_1^{h-1}D_2+a'(x_1)x_1^hD_2\big)-2if'(x_1)f(x_1)L,$$
and its localized operator $P^{(2)}_\varrho$ (see
\eqref{locgrushin}) at $\varrho=(0,x_2,0,\xi_2)$ is given by
$$(1+f(0)^2)D_t^2+a(0)^2t^{2h}\xi_2^2-(1-f(0)^2)ha(0)t^{h-1}\xi_2=$$
\begin{equation}\label{lockoh}
(D_{t}+i a(0)t^{h}\xi_{2})(D_{t}+i
a(0)t^{h}\xi_{2})^\ast+f(0)^2(D_{t}+i
a(0)t^{h}\xi_{2})^\ast(D_{t}+ia(0)t^{h}\xi_{2}).
\end{equation}. Hence, if
$f(0)\neq 0$, then $\textrm{Ker}\, P^{(2)}_\varrho=\textrm{Ker}\,
(D_{t}+i a(0)t^{h}\xi_{2})^\ast\cap\textrm{Ker}\, (D_{t}+i
a(0)t^{h}\xi_{2})=\langle 0\rangle$ so that $P^{(2)}_\varrho$ is
$L^2(\R)-$injective and $P$ turns out to be h.e., in $S_0$, with loss of $2h/(h+1)$
derivatives (see \eqref{locgrushin}). Otherwise, if we assume that
$f(0)=0$ and, precisely, that $0$ is a zero of order $k$ of $f$,
then $\textrm{Ker}\, P^{(2)}_\varrho=\textrm{Ker}\, (D_{t}+i
a(0)t^{h}\xi_{2})^\ast$ and a direct computation shows that
\begin{equation}\label{nucleo}
\big(D_{t}+i a(0)t^{h}\xi_{2}\big)^\ast\phi=0 \Longleftrightarrow
\phi=c e^{-\frac{a(0)}{h+1}t^{h+1}\xi_2}.
\end{equation}
Therefore, if $h$ is even, $P^{(2)}_\varrho$ is
$L^2(\R)-$injective and $P$ is again h.e., in $S_0$, with loss of $2h/(h+1)$
derivatives. On the other hand, if $h$ is odd, $P^{(2)}_\varrho$
has a non-trivial $L^2(\R)-$kernel, and hence $P$ can be
hypoelliptic only with  loss of derivatives larger than
$2h/(h+1)$. In order to apply Theorem \ref{teoprinc} for $h$ odd,
we need to compute the terms $\sum_{j\geq
0}\ell_{\frac{2}{h+1}-\frac{j}{h+1}}$ of Theorem \ref{main} and
defined by the iterative equations \eqref{psi'1}, \eqref{ell},
\eqref{psi'2}. Note that $P$ in \eqref{kext} does not depend on
the $x_2-$variable, so that the same holds for the terms
$\ell_{\frac{2}{h+1}-\frac{j}{h+1}}$.\\ Arguing as in the example
above, we see that the lowest eigenvalue $\l_{0}(\xi_2)$
identically vanishes, so that $\ell_{\frac{2}{h+1}}(\xi_2)\equiv
0$. Actually we are going to show that
\begin{equation}\label{claim}
\ell_{\frac{2}{h+1}-\frac{j}{h+1}}(\xi_2)\equiv 0\;\;\textrm{for
any}\; j<2k\;\;\textrm{and}\;\;
\ell_{\frac{2}{h+1}-\frac{2k}{h+1}}(\xi_2)\neq 0,
\end{equation}
where $k$ denotes the order of $x_1=0$ as zero of $f$.\\
Unfortunately, the direct computation of the terms
$\ell_{\frac{2}{h+1}}(\xi_2)$ by using iteratively the formulas
\eqref{psi'1}, \eqref{ell}, \eqref{psi'2} cannot be carried out;
indeed, if $h>1$, the eigenfunctions of the localized
operator \eqref{lockoh} seem to satisfy no iterative relations
similar to the classical ones verified by the standard Hermite
polynomials (see \cite{gun}). As a consequence, we are not able to
compute explicitly the action of
$E^{-2}_{(x_2,\xi_2)}=(P^{(2)}_\varrho)^{-1}$ and of the localized
operators $P_\varrho^{(2+j)}$ of order higher than $2$
(as done, for instance, in Section 6 \cite{Parenti3}). For that
reason, we follow a different approach. Firstly, by applying
Theorem \ref{main} to the operator $LL^*$ (i.e., the first
addendum in \eqref{kext}) we see that its hypoellipticity is
equivalent to the h.e. of an one-dimensional pseudodifferential
operator $\tilde{\ell}\sim\sum_{j\geq
0}\tilde{\ell}_{\frac{2}{h+1}-\frac{j}{h+1}}$. Secondly, since
$f(x_1)$ vanishes to the order $k$ at $x_1=0$, we observe that the
term $\big(f(x_1)L\big)^*f(x_1)L$ in \eqref{kext} gives no
contribution  in the computation of the
$\ell_{\frac{2}{h+1}-\frac{j}{h+1}}$ for $j<2k$ so that
$\ell_{\frac{2}{h+1}-\frac{j}{h+1}}=\tilde{\ell}_{\frac{2}{h+1}-\frac{j}{h+1}}$
for $j<2k$. Let us now proceed by contradiction. Assume that there
exists $0<j_0<2k$ such that
$\tilde{\ell}_{\frac{2}{h+1}-\frac{j}{h+1}}(\xi_2)\equiv 0$ if
$j<j_0$ and
$\tilde{\ell}_{\frac{2}{h+1}-\frac{j_0}{h+1}}(\xi_2)\neq 0$; thus,
$\tilde{\ell}(x_2,D_{x_2})$ would be an elliptic operator of order
$\frac{2}{h+1}-\frac{j_0}{h+1}$ and hence, trivially,
hypoelliptic. Therefore, from Theorem \ref{main} applied to
$LL^*$, it would follow that $LL^*$ is hypoelliptic (with loss of
$\frac{2h}{h+1}+\frac{j_0}{h+1}$ derivatives), but this is false.
To see this it is enough to construct a non smooth solution
$u(x_1,x_2)$ of $LL^*u(x_1,x_2)=0$. Since $a(0)>0$ and $h$ is odd, upon setting
$G(x_1)=\int_0^{x_1}g(t)dt=\int_0^{x_1}a(t)t^h dt$ we have that
$G(x_1)\geq c x_1^{h+1}\geq 0$ for $x_1$ near zero,  $c$ being a
suitable positive constant. The function
$$u(x_1,x_2)=\int_0^{+\infty}e^{ix_2\xi_2}e^{-G(x_1)\xi_2}\frac{d\xi_2}{1+\xi_2^4}$$
is well-defined near the origin, because  $e^{-G(x_1)\xi_2}\leq 1$, and solves $LL^*u(x_1,x_2)=0$ but it is not $C^\infty$ near the origin.
This proves the first part of \eqref{claim}. As a matter of fact, this also shows that $\tilde{\ell}_{\frac{2}{h+1}-\frac{2j}{h+1}}\equiv 0$ for every positive integer $j$.\\
We are left to check that
$\ell_{\frac{2}{h+1}-\frac{2k}{h+1}}(\xi_2)\neq 0.$ Since $f(x_1)$
vanishes to the order $k$ at $x_1=0$, we have that
$f(x_1)=b(x_1)x_1^k$ with $b\in C^w, b(0)\neq 0$; furthermore, a
careful analysis of the formulas \eqref{psi'1}, \eqref{ell},
\eqref{psi'2} shows that
$\ell_{\frac{2}{h+1}-\frac{2k}{h+1}}(\xi_2)$ is the sum of
$\tilde{\ell}_{\frac{2}{h+1}-\frac{2k}{h+1}}(\xi_2)\equiv 0$ and
of the contribution due to the operator
$\big(f(x_1)L\big)^*f(x_1)L$. Hence, recalling \eqref{nucleo}, we
get that
$$\ell_{\frac{2}{h+1}-\frac{2k}{h+1}}(\xi_2)=\<b(0)x_1^k\big(D_{t}+i a(0)t^{h}\xi_{2}\big)^\ast \big(b(0)x_1^k(D_{t}+i a(0)t^{h}\xi_{2}\big)\phi,\phi\>$$
$$=\|b(0)x_1^k\big(D_{t}+i a(0)t^{h}\xi_{2}\big)\phi\|_0^2>0.$$
This completes the proof of our claim \eqref{claim}. Therefore
$\ell(x_2,D_{x_2})$ turns out to be an elliptic operator of order
$\frac{2}{h+1}-\frac{2k}{h+1}$ and hence, as operator in
$\textrm{OP}S^{\frac{2}{h+1}}$, is hypoelliptic with loss of
$\frac{2k}{h+1}$. As a consequence of Theorem \ref{main}, $P$ is
hypoelliptic with loss of $\frac{2(h+k)}{h+1}$ in $S_0$.\\
Arguing as above near any zero of $g(x_1)$ we obtain the following
proposition.\\
\begin{proposition}\label{hkext}
Assume that $g$ is a polynomial or a analytic function with a finite number of  real zeroes. Set $R_g=\{x_1\in\R\,|\, g(x_1)=0\}=\{r_1,r_2,...,r_k\}$ and denote by $h_j^g\; (j=1,2,...,k)$ the order $r_j$ as zero of $g$. Moreover, define the quantities $h_j^f$
as follows. Set $h_j^f=0$ if $r_j$ is a zero of $g$ of even order or if $f(r_j)\neq 0$. Otherwise, set $h_j^f$ to be the order of $r_j$ as zero of $f$. Then $P$ is
$C^\infty-$hypoelliptic with loss of derivatives given by
$$\max_{j=1,...,k}\left\{\frac{2h_j^g+2h_j^f}{h_j^g+1}\right\}.$$
\end{proposition}
We point out that for a general real analytic function $g$, the set $R_g$ is not necessarily finite. Via the Weierstrass theorem we can construct two analytic functions $g,f$ having infinite real zeroes of arbitrary multiplicities, so that the quantity
$\sup\left\{\frac{2h_j^g+2h_j^f}{h_j^g+1}\right\}$ can be finite or infinite depending on such zeroes. However also in this case the operator $P$ is h.e.
in any bounded open set.

\subsection{A sum of squares of complex vector fields:}\label{uniquefield}
We complete this section by analyzing the following sum of squares of a complex and a real vector field in $\R^2$:
\begin{equation}\label{squares}
P=\big(D_{x_1}-i x_1^{k_1} D_{x_2}\big)^\ast\big(D_{x_1}-i x_1^{k_1} D_{x_2}\big)+\big(x_1^{k_2} D_{x_2})^\ast \big(x_1^{k_2} D_{x_2}),\quad k_1,k_2\in\N.
\end{equation}
It is worth noting that the operator $ P $ satisfies the complex H\"ormander condition, i.e. the brackets of the fields $D_{x_1}-i x_1^{k_1} D_{x_2}$ and $x_1^{k_2} D_{x_2}$
of length up to $ k_2 + 1 $ generate a two dimensional complex Lie algebra. Obviously, also the real H\"ormander condition (using as vector fields both the real and the imaginary parts of the vector fields
defining $ P $)  is satisfied. We are going to show that also the presence of a single genuine complex vector field can have a strong impact on the loss of
regularity of an operator in a sum of squares form. Indeed, a sum of squares of real vector fields (satisfying the H\"ormander condition) is actually subellitic
(i.e., the loss of derivatives is always less than $2$, see Theorem 22.2.1 \cite{ho}); whereas, in our model, the loss of regularity can be arbitrarily large if
$k_1$ is odd and $k_2>>k_1$. This phenomenon was firstly pointed out by
J.J.Kohn in \cite{kohn} in the case of two complex vector fields in $\R^2$.\\
We collect our results in the proposition below.
\begin{proposition}\label{propunique} We have the following statements:\\
$(i)$ if  $0<k_2\leq k_1$, then $P$ is hypoelliptic with loss of $2k_2/(k_2+1)$ derivatives;\\
$(ii)$ if $k_2>k_1$ and $k_1$ is an even integer, $P$ is hypoelliptic with loss of $2k_1/(k_1+1)$ derivatives;\\
$(iii)$ if $k_2>k_1$ and $k_1$ is an odd integer, $P$ is hypoelliptic with loss of $2k_2/(k_1+1)$ derivatives.
\end{proposition}
{\bf Proof.} If  $0<k_2\leq k_1$, then $P$ can be written in the form \eqref{grushin} with $h=k_2$ and the corresponding localized operator at $\v=(x_2,\xi_2)$ is
given by $P_\v=D_t^2+t^{2k_2}\xi_2^2$ if $k_2<k_1$ and by $P_\v=\big(D_{t}-i t^{k_1} \xi_2\big)^\ast\big(D_{t}-i t^{k_1} \xi_2\big)+t^{2k_2}\xi_2^2$ if $k_1=k_2$.
By taking the $L^2(\R_t)$-scalar $\langle P_\v u,u\rangle$ it is easily seen that $P_\v$ is injective in $L^2(\R_t)$ so that $P$ is hypoelliptic with loss of
$2k_2/(k_2+1)$ derivatives (see, for instance, Thm. 1.12 \cite{marco}).\\
If $k_2>k_1$, by choosing $h=k_1$ the localized operator at $\v=(x_2,\xi_2)$ is given by
$$P_\v=\big(D_{t}-i t^{k_1} \xi_2\big)^\ast\big(D_{t}-i t^{k_1} \xi_2\big),$$
and is injective in $L^2(\R_t)$ for even $k_1$, as seen in \eqref{nucleo}. Thus the statement $(ii)$ readily follows.\\
Finally, if $k_1$ is an odd integer, we see that, for $\xi_2>0$,
$$\phi=e^{-\frac{t^{k_1+1}}{k_1+1}\xi_2}\in \textrm{Ker}_{L^2(\R_t)} P_\v,$$
and, since $P_\v^\ast=P_\v$, from Remark \ref{selfadj} we gets that $\phi_1=\phi_2=\phi$. Moreover, due to the structure of $P$, we see that all the high order
localized operators $P_\v^{(2+r)}$ (see \eqref{lochigh}) identically vanish, except for $r=2(k_2-k_1)$ for which we have $P_\v^{(2+2(k_2-k_1))}=t^{2k_2}\xi_2^2$. By
applying Theorem \ref{main}, we are reduce to compute the symbol $\sum_{j\geq 0}\ell_{\frac{2}{k_1+1}-\frac{j}{k_1+1}}$. Note that $P$ in \eqref{squares} does not
depend on the $x_2-$variable, so that the same must hold for the quantities involved in the terms $\ell_{\frac{2}{h+1}-\frac{j}{h+1}}$. An inspection of the
equation \eqref{ell} shows that
\begin{equation*}
\ell_{\frac{2}{k_1+1}-\frac{j}{k_1+1}}\equiv 0\;\;\;\textrm{for}\;\; 0\leq j<2(k_2-k_1)\quad\textrm{and}\quad
\ell_{\frac{2}{k_1+1}-\frac{2(k_2-k_1)}{k_1+1}}=\xi_2^2\|t^{k_2}\phi\|^2_0\neq 0.
\end{equation*}
As a consequence, $\ell(x_2,D_2)$ is hypoelliptic with loss of $\frac{2(k_2-k_1)}{k_1+1}$ derivatives, and, in view of Theorem \ref{main}, this concludes the proof.
\subsection{An example with a large loss of regularity in the even case}
In Section \ref{koext} and \ref{uniquefield} we have shown operators which can be h.e. with an arbitrary large loss of derivatives provided that $h$ in
\eqref{grushin} be an odd integer. This is not a true restriction and it is easy to recognize the same behavior in the following example:
\begin{equation}\label{gteven}
D_{x_1}^2+x_1^{2h}D_{x_2}^2+\big(x_1^k+h+1\big)x_1^{h-1}D_{x_2}
\end{equation}
where $h$ is even and $k\in\N$ is odd. Again, by Theorem \ref{main} it is enough to compute the symbol $\sum_{j\geq 0}\ell_{\frac{2}{h+1}-\frac{j}{h+1}}$, similarly
as done in the previous section. Note that $\ell_{\frac{2}{h+1}}\equiv 0$ due to \eqref{avreven} with $j_0=0$ and that
\begin{equation*}
\ell_{\frac{2}{h+1}-\frac{j}{h+1}}\equiv 0\;\;\;\textrm{for}\;\; 0\leq j<k\quad\textrm{and}\quad
\ell_{\frac{2}{h+1}-\frac{k}{h+1}}=|\xi_2|\cdot\|t^{(k+h-1)/2}\phi\|^2_0\neq 0.
\end{equation*}
Therefore \eqref{gteven} is h.e. with loss of $\frac{2h+k}{h+1}$ derivatives.

\section{Tangential Grushin type operators}\label{mod}
In this section we consider a special type of operator
\eqref{grushin}, whose coefficients depend only by the tangential
variables to $\textrm{Char }P$, i.e.
\begin{equation}\label{tangrushin}
 P=D_{x_1}^2+a(x'){x_1}^{2h}\Delta_{x'}+{x_1}^{h-1}B(x',D_{x'}),\qquad x_1\in\R,\;x'=(x_2,...,x_{n})\in\R^{n-1}.
\end{equation}
Throughout this section we call $P$ the tangential Grushin operator.\\
We assume that the classical condition of $C^\infty$-hypoellipticity with  loss of $2h/(h+1)$ derivatives is violated at a point $(x_0',\xi_0')$; namely, there is
an eigenvalue $\l_{j_0}$ of the localized operator $P^{(2)}_\varrho$ that vanishes at $\varrho=(x_0',\xi_0')$ (see \eqref{avr}). We shall prove that the h.e. of $P$
is strictly related to the sign of $1/i\{\l_{j_0},\overline{\l_{j_0}}\}(x_0',\xi_0')$, similarly to what happens in the isotropic case $h=1$ (see Thm. 1.2
\cite{h2}). Here the main difficulty concerns the fact that, unlike the case $h=1$,
no explicit formula is known for the  eigenvalue $\l_{j_0}$.\\
However, it is worth noting that the h.e. (with minimal loss) of a general Grushin-type operator \eqref{grushin} (i.e. not tangential) does not depend on the
Poisson brackets condition \eqref{autopoisson} for $h\neq 3$ (as shown in Theorem \ref{teoprinc}); furthermore, also for $h=3$, that condition is no longer
necessary. For instance, if we choose $k_1=3$ odd and $k_2=4$ the operator \eqref{squares}  turns out to be h.e. with loss of $2k_1/(k_1+1)+1/2=2$
although in this case we have that $\{\l_{j_0},\overline{\l_{j_0}}\}\equiv 0$.\\
Let us apply Theorem \ref{main} to the operator \eqref{tangrushin} and briefly discuss the
structure of the operator $\ell(x',D_{x'})$, with $\ell(x',\xi')\sim \sum_{j\geq 0}\ell_{\frac{2}{h+1}-\frac{j}{h+1}}(x',\xi')$ associated with $P$ via
Theorem \ref{main} and defined by the iterative formulas \eqref{psi'1},
\eqref{ell}, \eqref{psi'2}.\\
Firstly, let us show that
\begin{equation}\label{lvanish}
\ell_{\frac{2}{h+1}-\frac{j}{h+1}}(x',\xi')\equiv 0,\quad\textrm{if}\quad 0<j\leq h.
\end{equation}
In view of \eqref{ell}, we have that
$$
\ell_{2/(h+1)-j/(h+1)}(x',\xi')=\Bigl(P_{(x',\xi')}^{(2)}{\psi'}_{-j/(h+1),1} (x',\xi';\cdot),\phi_2(x',\xi';\cdot)\Bigr)_{L^2(\R)}+$$
$$-\sum_{\text{\scriptsize$\begin{matrix}(h+1)|\alpha|+q=j\\ 0\leq q<j\end{matrix}$}}\hspace{-.2cm}
\frac{1}{\alpha !i^{|\alpha|}}\p_{x'}^\alpha\ell_{2/(h+1)-q/(h+1)}(x',\xi')\Bigl(
\p_{\xi'}^\alpha\phi_{2}(x',\xi';\cdot),\phi_2(x',\xi';\cdot)\Bigr)_{L^2(\R)}\hspace{-.3cm}+$$
\begin{equation}\label{ell2}
+\sum_{\text{\scriptsize$\begin{matrix}(h+1)|\alpha|+s+q=j\\
0\leq q<j\end{matrix}$}} \frac{1}{\alpha !i^{|\alpha|}}\Bigl(\p_{\xi'}^\alpha P_{(x',\xi')}^{(2+s)}\p_{x'}^\alpha{\psi'}_{-q/(h+1)}
(x',\xi';\cdot),\phi_2(x',\xi';\cdot)\Bigr)_{L^2(\R)}.
\end{equation}
Note that the first and the second term in the r.h.s. are identically zero  if $0<j\leq h$; indeed, from \eqref{psi'1} it follows that
$$
{\psi'}_{-j/(h+1),1}= -\sum_{\text{\scriptsize$\begin{matrix}(h+1)|\alpha|+q=j\\ 0\leq q<j\end{matrix}$}} \frac{1}{\alpha
!i^{|\alpha|}}\Bigl(\p_{x'}^\alpha{\psi'}_{-q/(h+1)}, \p_{\xi'}^\alpha{\phi}_{1}\Bigr)_{L^2(\R_{x_1})}\cdot{\phi}_{1};
$$
the condition $(h+1)|\alpha|+q=j\leq h$ trivially yields $\a=0$ and $q=j$, so that there are no terms in the above sum for $0\leq q<j$ and hence
${\psi'}_{-j/(h+1),1}\equiv 0$.\\
The same argument applies to the second term in the r.h.s. of \eqref{ell2}. Finally, as for the last term in \eqref{ell2}, note that all the localized operators
$P_{(x',\xi')}^{(2+s)}$ of the tangential Grushin operator are identically zero for any $s>0$ (see \eqref{lochigh}); therefore, from
$(h+1)|\alpha|+s+q=j\leq h$ it follows that $s=0, \a=0$ and hence, due to the condition $q<j$, there are no terms in the last sum of \eqref{ell2}.\\
Therefore, in view of \eqref{lvanish}, we can use Proposition \ref{tangential} (in the appendix I) in order to study the hypoellipticity of $\ell(x',D_{x'})$ (and so the h.e. of $P$). To this purpose, we are left to compute its principal symbol $\ell_{\frac{2}{h+1}}(x',\xi')=
\<P^{(2)}_{(x',\xi')}\phi_1,\phi_2\>_{L^2(\R)}$, where $\phi_1, \phi_2$ are defined in Lemma \ref{quad}. Recalling Lemma \ref{quad} and Remark \ref{selfadj2} we
have
\begin{eqnarray*}
\ell_{\frac{2}{h+1}}(x',\xi')&=&\<P^{(2)}_{(x',\xi')}\phi_1,\phi_2\>_{L^2(\R)}=\<P^{(2)}_{(x',\xi')}\big(\phi-(\phi-\phi_1)\big),\phi_2\>_{L^2(\R)}\\
&=&\<P^{(2)}_{(x',\xi')}\phi,\phi_2\>_{L^2(\R)}+\<P^{(2)}_{(x',\xi')}\big(\phi-\phi_1\big),\phi_2\>_{L^2(\R)}\\
&=&\l_{j_0}(x',\xi')\<\phi,\phi_2\>_{L^2(\R)}+\<\phi-\phi_1,P^{(2)\ast}_{(x',\xi')}\phi_2\>_{L^2(\R)}\\
&=&\l_{j_0}(x',\xi')\<\phi,\phi_2\>_{L^2(\R)}+O(\|(x'-x_0',\xi'-\xi'_0)\|^2)
\end{eqnarray*}
since $\phi-\phi_1\rightarrow 0$ and $P^{(2)\ast}_{(x',\xi')}\phi_2\rightarrow \l_{j_0}(x_0',\xi_0')\phi(x_0',\xi_0')=0$ as soon as $(x',\xi')\rightarrow
(x_0',\xi_0')$.\\
In order to apply Prop. \ref{tangential}, we need to compute the $T^\ast\R^{n-1}-$Poisson brackets
$\frac{1}{i}\{\ell_{\frac{2}{h+1}},\overline{\ell_{\frac{2}{h+1}}}\}$ at $(x_0',\xi_0')$. Since $\l_{j_0}(x_0',\xi_0')=0$ and
$\phi\big((x_0',\xi_0'),\cdot\big)=\phi_2\big((x_0',\xi_0'),\cdot\big)$, we have
$$
\frac{1}{i}\{\ell_{\frac{2}{h+1}},\overline{\ell_{\frac{2}{h+1}}}\}(x_0',\xi_0')=
\frac{1}{i}\{\l_{j_0},\overline{\l_{j_0}}\}(x_0',\xi_0')\cdot
|\<\phi\big((x_0',\xi_0'),\cdot\big),\phi_2\big((x_0',\xi_0'),\cdot\big)\>_{L^2(\R)}|^2.
$$
$$=\frac{1}{i}\{\l_{j_0},\overline{\l_{j_0}}\}(x_0',\xi_0')\cdot
\underbrace{\|\phi\big((x_0',\xi_0'),\cdot\big)\|^4_{L^2(\R)}}_{> 0},$$ hence, in view of Prop. \ref{tangential}, the analysis of the hypoellipticity of $P$ is
reduced to the study the sign of $\frac{1}{i}\{\l_{j_0},\overline{\l_{j_0}}\}(x_0',\xi_0')$.\\
 Unlike the case $h=1$ (see \cite{h2}), here the eigenvalue $\l_{j_0}(x',\xi')$ cannot be explicitly computed since the
spectrum of the anharmonic oscillator is not known (see, for instance, \cite{marco2}). In order to overcome this difficulty, we use classical perturbation theory
(see \cite{kato} and \cite{RS}) to get a
convenient approximation of $\l_{j_0}(x',\xi')$ near $(x_0',\xi_0')$, enough to compute the Poisson brackets at $(x_0',\xi_0')$.\\
To this aim, let us write the localized operator $P_\varrho$ at  (see \eqref{locgrushin}) as a perturbation of $P_{\v_0}$ for $\v=(0,x',0,\xi')$ near
$\v_0=(0,x_0',0,\xi_0')$
$$P_{\varrho}=\underbrace{D_t^2+a(x_0')|\xi_0'|^2 t^{2h}+b_1(x_0',\xi_0')t^{h-1}}_{=\; \bold{ P_{\v_0}}}$$
\begin{equation}\label{locper}
+\underbrace{\big(a(x')|\xi'|^2-a(x_0')|\xi_0'|^2\big)}_{=\; \bold{A(\v)}}\cdot t^{2h}+ \underbrace{\big(b_1(x',\xi')-b_1(x_0',\xi_0')\big)}_{=\; \bold{B(\v)}}\cdot
t^{h-1}.
\end{equation}
We point out that $A(\v), B(\v)$ are smooth functions such that $A(\v_0)=0$ and $B(\v_0)=0$. It is worth noting that $A(\v)$ is a real smooth function.\\
We get the following lemma.
\begin{lemma}\label{conto}
One has
$$\l_{j_0}(x',\xi')=A(x',\xi')\langle t^{2h}\phi_{\v_0},\phi_{\v_0}\rangle_{L^2(\R_t)}+
B(x',\xi')\langle t^{h-1}\phi_{\v_0},\phi_{\v_0}\rangle_{L^2(\R_t)}$$
$$+ O\big(\|x'-x_0'\|^2+\|\xi'-\xi'_0\|^2\big)$$ where
$\phi_{\v_0}=\phi\big((x_0',\xi_0'),\cdot\big)\in\textrm{Ker}\ P_{\v_0}$ with $\|\phi_{\v_0}\|_{L^2(\R_t)}=1$ (see Remark \ref{selfadj2}).
\end{lemma}
{\bf Proof.} Since $\|\phi_{\v_0}\|^2\neq 0$, by continuity we have that $\langle \phi_{\v}, \phi_{\v_0} \rangle\neq 0$ if $\v$ is near $\v_0$. Moreover, from
\eqref{locper} and Remark \ref{selfadj2} we get
\begin{eqnarray*}
\l(\v)&=&\frac{\langle P_\v \phi_\v, \phi_{\v_0}\rangle}{\langle \phi_\v, \phi_{\v_0}\rangle}=\frac{\langle\big(P_{\v_0}+A(\v)t^{2h}+B(\v)t^{h-1}\big)\phi_\v,
\phi_{\v_0}\rangle}{\langle \phi_\v, \phi_{\v_0}\rangle}\\
&=&\frac{\langle\phi_\v, P_{\v_0}\phi_{\v_0}\rangle}{\langle \phi_\v, \phi_{\v_0}\rangle}+A(\v)\frac{\langle t^{2h}\phi_\v, \phi_{\v_0}\rangle}{\langle \phi_\v,
\phi_{\v_0}\rangle}+B(\v)\frac{\langle t^{h-1}\phi_\v, \phi_{\v_0}\rangle}{\langle \phi_\v, \phi_{\v_0}\rangle}\\
&=&A(\v)\frac{\langle t^{2h}\phi_{\v}, \phi_{\v_0}\rangle}{\langle \phi_\v, \phi_{\v_0}\rangle}+B(\v)\frac{\langle t^{h-1}\phi_\v, \phi_{\v_0}\rangle}{\langle
\phi_\v, \phi_{\v_0}\rangle}
\end{eqnarray*}
We have
$$A(\v)\frac{\langle t^{2h}\phi_{\v}, \phi_{\v_0}\rangle}{\langle \phi_\v, \phi_{\v_0}\rangle}=A(\v)\frac{\langle t^{2h}\phi_{\v_0}, \phi_{\v_0}\rangle+\langle t^{2h}\big(\phi_{\v}-\phi_{\v_0}\big), \phi_{\v_0}\rangle}
{1+\langle \big(\phi_{\v}-\phi_{\v_0}\big), \phi_{\v_0}\rangle}$$
$$=A(\v)\langle t^{2h}\phi_{\v_0},\phi_{\v_0}\rangle+ O\big(\|\v-\v_0\|^2\big),$$
and similarly
$$B(\v)\frac{\langle t^{h-1}\phi_\v, \phi_{\v_0}\rangle}{\langle
\phi_\v, \phi_{\v_0}\rangle}=B(\v)\langle t^{h-1}\phi_{\v_0},\phi_{\v_0}\rangle+ O\big(\|\v-\v_0\|^2\big),$$ which completes the proof.\fine
%The spectrum of $P_\v$ is discrete made of simple eigenvalues; furthermore, in view of \eqref{avr}, $\l_{j_0}(\v_0)=0$ is a simple isolated eigenvalue of $P_{\v_0}$. From Theorems XII.5, XII.6 and XII.8 \cite{RS} there exists a $ \delta> 0 $, small enough, such that the operator
%\begin{equation}\label{formulaint}\Pi_{\v} = \frac{1}{2\pi i} \oint_{|\sigma | = \delta} (\sigma I - P_\v)^{-1} d\sigma\end{equation}
%is the orthogonal projection onto the eigenspace $\textrm{Ker} \big(P_\v -\lambda_{j_0}(\v )I\big)$, provided that $\v$ is sufficiently near $\v_0$.

\noindent Taking into account that $A(\v)$ is a real function, by the above lemma we obtain
\begin{align}\label{Poissonesp}
\frac{1}{i}\{\l_{j_0},\overline{\l_{j_0}}\}(x_0',\xi_0')&=2\langle t^{2h}\phi_{\v_0},\phi_{\v_0}\rangle \langle t^{h-1}\phi_{\v_0},\phi_{\v_0}\rangle \textrm{Im}
\{B,A\}(x_0',\xi_0')\nonumber\\
&+\frac{1}{i}\langle t^{h-1}\phi_{\v_0},\phi_{\v_0}\rangle^2 \{B,\overline{B}\}(x_0',\xi_0')\nonumber\\
&=2\langle t^{2h}\phi_{\v_0},\phi_{\v_0}\rangle \langle t^{h-1}\phi_{\v_0},\phi_{\v_0}\rangle \textrm{Im} \big\{b_1(x',\xi'),a(x')|\xi'|^2\big\}(x_0',\xi_0')\nonumber\\
&  +\big(\langle t^{h-1}\phi_{\v_0},\phi_{\v_0}\rangle\big)^2 \frac{1}{i}\big\{b_1(x',\xi'),\overline{b_1(x',\xi')}\big\}(x_0',\xi_0').
\end{align}
In order to study its sign, we do not need to know explicitly the eigenfunction $\phi_{\v_0}$; it is enough to compute $\langle
t^{2h}\phi_{\v_0},\phi_{\v_0}\rangle$ in terms of $\langle t^{h-1}\phi_{\v_0},\phi_{\v_0}\rangle$. To this aim, we exploit a classical scaling argument. If we
perform the change of variable $t=sy$ ($s\in\R$) in the localized operator $P_{\v_0}=D_t^2+a|\xi_0'|^2 t^{2h}+b_1 t^{h-1}$, we get
$$P_{\v_0,s}=\frac{1}{s^2}D_y^2+a|\xi_0'|^2 s^{2h}y^{2h}+b_1 s^{h-1}y^{h-1},\quad \phi_{\v_0, s}(y):=\phi_{\v_0}(sy).$$
As a consequence of the spectral invariance, we have the identity $P_{\v_0,s}\phi_{\v_0, s}=0$ and a differentiation with respect to $s$ yields
\begin{equation*}
\big(-\frac{2}{s^3}D_y^2+a|\xi_0'|^2 2h s^{2h-1}y^{2h} +b_1(h-1) s^{h-2}y^{h-1}\big)\phi_{\v_0, s}+P_{\v_0,s}(\partial_s\phi_{\v_0, s})=0.
\end{equation*}
Taking the scalar product with $\phi_{\v_0, s}$ and choosing $s=1$ give
$$-2\|D_t\phi_{\v_0}\|_{L^2(\R_t)}^2+2h a|\xi_0'|^2\langle t^{2h}\phi_{\v_0},\phi_{\v_0}\rangle+(h-1)b_1\langle  t^{h-1}\phi_{\v_0},\phi_{\v_0}\rangle$$
\begin{equation}\label{scaling}
+\langle {\partial_s\phi_{\v_0,s}}_{|s=1}, \underbrace{P^*_{\v_0}(\phi_{\v_0})}_{{{=P_{\v_0}(\phi_{\v_0}) = 0}}}\rangle=0.
\end{equation}
Here, for the sake of simplicity, we have replaced $y$ by
$t$ since $t=y$ for $s=1$. By taking  the scalar product with $\phi_{\v_0}$ in the identity $P_{\v_0}\phi_{\v_0}=0$, we obtain
$$\|D_t\phi_{\v_0}\|_{L^2(\R_t)}^2+ a|\xi_0'|^2\langle t^{2h}\phi_{\v_0},\phi_{\v_0}\rangle+b_1\langle  t^{h-1}\phi_{\v_0},\phi_{\v_0}\rangle=0,$$
whence, together with \eqref{scaling} we finally get
$$\langle t^{2h}\phi_{\v_0},\phi_{\v_0}\rangle=-\frac{b_1(x',\xi')}{2a(x')|\xi_0'|^2}\langle t^{h-1}\phi_{\v_0},\phi_{\v_0}\rangle.$$
As a trivial consequence, note that $\langle t^{h-1}\phi_{\v_0},\phi_{\v_0}\rangle\neq 0$ since $\langle
t^{2h}\phi_{\v_0},\phi_{\v_0}\rangle=\|t^h\phi_{\v_0}\|_0^2\neq 0$. By using the above formula in \eqref{Poissonesp} we see that
$$\frac{1}{i}\{\l_{j_0},\overline{\l_{j_0}}\}(x_0',\xi_0')=\big(\langle t^{h-1}\phi_{\v_0},\phi_{\v_0}\rangle\big)^2\cdot$$
$$\Big(\frac{1}{i}\big\{b_1(x',\xi'),\overline{b_1(x',\xi')}\big\}(x_0',\xi_0')-\frac{b_1(x_0',\xi_0')}{a(x_0')|\xi_0'|^2}
 \big\{\textrm{Im}\, b_1(x',\xi'),a(x')|\xi'|^2\big\}(x_0',\xi_0')\Big).$$
Therefore, from this identity and Proposition \ref{tangential} the following result is proved.
\begin{proposition}\label{proptan}
Let P be as in \eqref{tangrushin} such that \eqref{avreven}, \eqref{avrodd} are satisfied, according to the parity of $h$. Then $P$ can only be hypoelliptic with
loss $\sigma \geq \frac{2h}{h+1}+\frac{1}{2}$. This threshold is attained if and only if
\begin{equation}\label{poissonesplicito}
\frac{1}{i}\big\{b_1(x',\xi'),\overline{b_1(x',\xi')}\big\}(x_0',\xi_0')- \frac{b_1(x_0',\xi_0')}{a(x_0')|\xi_0'|^2}
 \big\{\mathrm{Im}\, b_1(x',\xi'),a(x')|\xi'|^2\big\}(x_0',\xi_0')<0
\end{equation}
\end{proposition}
Let us complete this section by pointing out a peculiarity of the even case if $P$ is a differential operator (i.e., $b_1(x_1,x',\xi')$ is homogeneous in $\xi'$).
Assume that $h$ be an even integer; moreover, set $\v=(x',\xi')$ and, accordingly, $-\v=(x',-\xi')$. Due to the parity of $h$, the operators $P_{-\v}$ and $P_{\v}$
are unitarily equivalent, via the change of variable $x\rightarrow -x$. As a consequence, they have the same spectrum and from \eqref{ord} one easily gets that
\begin{equation}\label{evenp}
\l_{j_0}(x',-\xi')=\l_{j_0}(x',\xi')\qquad \textrm{for every }\quad (x',\xi')\in T^\ast\R^{n-1}\setminus 0.
\end{equation}
Furthermore, in view of \eqref{avreven}, one has that $\l_{j_0}(x'_0,-\xi'_0)=0=\l_{j_0}(x'_0,\xi'_0)$. Hence, if we had
$\frac{1}{i}\{\l_{j_0},\overline{\l_{j_0}}\}(x_0',\xi_0')<0$, then from \eqref{evenp} it would follow that
$\frac{1}{i}\{\l_{j_0},\overline{\l_{j_0}}\}(x_0',-\xi_0')>0$. Therefore, we have the following remark.
\begin{remark}
If $h$ is an even integer and $P$ in \eqref{tangrushin} is a differential operator, then $P$ can never be hypoelliptic with loss of
$\sigma=\frac{2h}{h+1}+\frac{1}{2}$ derivatives. However, $P$ can be microhypoelliptic at $\v_0=(0,x'_0,0,\xi'_0)$ (with such a loss of regularity) if and only if
\eqref{poissonesplicito} holds.
\end{remark}
Let us now apply Prop. \ref{proptan} to the following example in $\R^3$:
$$P=D_1^2+x_1^{2h}\big(D_2^2+D_3^2)+\a x_1^{h-1}D_2+f(x_2,x_3)x_1^{h-1}D_3,$$
where $\a\in\mathbb{C}: \mathrm{Im}\ \a\neq 0, f\in \textrm{C}^\infty(\R^2,\R)$. In view of the above remark, let us assume that $h$ is an odd integer and suppose
that \eqref{avrodd} holds at $\v_0=(\tilde{x}_2,\tilde{x}_3,\tilde{\xi}_2,\tilde{\xi}_3)$, i.e. for some positive integer $j_0$
$$\a\tilde{\xi}_2+f(\tilde{x}_2,\tilde{x}_3)\tilde{\xi}_3=(-1)^{j_0}|(\tilde{\xi}_2,\tilde{\xi}_3)|-(h+1)(j_0+\theta(j_0))
|(\tilde{\xi}_2,\tilde{\xi}_3)|.$$ Since the r.h.s. is real, we have that $\mathrm{Im} \a\ \tilde{\xi}_2=0 \Longrightarrow \tilde{\xi}_2=0$, whence we
get
$$f(\tilde{x}_2,\tilde{x}_3)=(-1)^{j_0}-(h+1)(j_0+\theta(j_0)),\quad (\tilde{\xi}_2,\tilde{\xi}_3)=(0,\tilde{\xi}_3)$$
where we have assumed, without loss of generality, that $\tilde{\xi}_3>0$. Finally, from Prop. \ref{proptan} it follows that $P$ is hypoelliptic with loss
$\frac{2h}{h+1}+\frac{1}{2}$ if and only if
$$\mathrm{Im} \a\cdot \frac{\partial f}{\partial x_2}(\tilde{x}_2,\tilde{x}_3)<0.$$

\setcounter{section}{0}
\renewcommand\thesection{\Alph{section}}
\section{Appendix I.  Hypoellipticity for $h$-homogeneous symbols}\label{simbani}
\setcounter{equation}{0}
\setcounter{theorem}{0}

\renewcommand{\thetheorem}{\thesection.\arabic{theorem}}
\renewcommand{\theproposition}{\thesection.\arabic{proposition}}
\renewcommand{\thelemma}{\thesection.\arabic{lemma}}
\renewcommand{\thedefinition}{\thesection.\arabic{definition}}
\renewcommand{\thecorollary}{\thesection.\arabic{corollary}}
\renewcommand{\theequation}{\thesection.\arabic{equation}}

In this appendix we carry out the study of the
$C^\infty-$hypoellipticity of operators
modelled on $L$ in the statement $4$ of Theorem \ref{main}.\\
To this aim, let $a(y,\eta)\in
S_{1,0}^{m'}(\R_y^\nu\times\R_\eta^\nu)$ be a symbol with the
(poly)homogenous asymptotic expansion (see Def.18.1.5 Vol.III \cite{ho})
\begin{equation}\label{amj}
a\sim\sum_{j\geq 0}a_{m'-\frac{j}{h+1}},\qquad
a_{m'-\frac{j}{h+1}}(y,\l\eta)=\l^{m'-\frac{j}{h+1}}a_{m'-\frac{j}{h+1}}(y,\eta)
\;\; \;\textrm{for any}\; \l>1,
\end{equation}
and $A=a(y,D_y)$ the corresponding properly supported operator.\\
From now on we assume that $h>1$, for the case $h=1$ is treated in
\cite{h2}. We can rephrase the hypoellipticity of $A=a(y,D_y)$ in
terms of suitable apriori estimates, as shown by the following
classical result.
\begin{lemma}\label{ine}
If $A$ is hypoelliptic with loss of $\sigma\geq 0$ derivatives,
then for any $s\in\R, \mu\in\R: \mu<s+m'-\sigma$ and any compact
set $K\subset\R^\nu$ there exists a positive constant
$C=C(K,s,\mu)$ such that
\begin{equation}\label{stimaauto}
\|u\|_{s+m'-\sigma}\leq C\big(\|Au\|_s+\|u\|_\mu\big),\qquad
\forall u\in C^\infty_0(K).
\end{equation}
Furthermore, assume now that $0\leq\sigma<1$. If $A$ satisfies the
inequality \eqref{stimaauto} for $s=0$ and any compact set $K$,
then $A$ is hypoelliptic with  loss of $\sigma$ derivatives.
\end{lemma}
{\bf Proof.} Since $A$ is properly supported, given any compact
set $K\subset \R^\nu$ there exists a compact set $\tilde{K}\subset
\R^\nu$
 such that $\mr{supp}(Au)\subseteq \tilde{K}$ if $u \in
\cal{E'}(\mathrm{K})$. Fix $s,\mu \in \R$ with $\mu<s+m'-\sigma$
and define
$$H_{s,\mu}(K)=\{u \in H^{\mu}\cap {\cal{E'}(\mathrm{K})}\; | \; Pu \in H^s\cap {\cal{E'}(\tilde{\mathrm{K}})} \}$$
equipped with the norm $[u]_{s,\mu}=\|Pu\|_s+\|u\|_{\mu}$. In view
of the hypoellipticity of $A$, the space $H_{s,\mu}(K)$ is embedded in
$H^{s+m'-\sigma}(K)$ and the imbedding is closed;
therefore an application of the Closed Graph Theorem yields inequality \eqref{stimaauto}.\\
The converse of the statement requires a sharp regularization argument and is a straightforward consequence of Lemma 27.1.5 in Vol.IV \cite{ho}.\\
The following lemma shows which is the minimal loss of derivatives
expected from $A$ whenever $A$ is not elliptic, i.e. ${\rm
Char}(A)=\{(y,\eta)\in T^\ast\R^\nu\setminus 0,
a_{m'}(y,\eta)=0\}\neq\emptyset$.
\begin{lemma}\label{loss}
Suppose that $a_{m'}(y_0,\eta_0)=0$ for some $(y_0,\eta_0)\in
T^\ast\R^\nu\setminus 0$, then $A$ can be hypoelliptic only with
loss $\sigma\geq\frac{1}{h+1}$ of derivatives. If, furthermore, $h$ is an
odd integer and  $a_{m'-\frac{1}{h+1}}(y_0,\eta_0)=0$,  then $\sigma\geq\frac{2}{h+1}$.
\end{lemma}
{\bf Proof.} If $A$ is hypoelliptic but not elliptic, it must lose
derivates. Assume that $A$ is hypoelliptic with loss of $\sigma$
derivatives, then by Lemma \ref{ine} ones has that
\begin{equation}\label{pal}
\|u\|_{s+m'-\sigma}\leq C\big(\|Au\|_s+\|u\|_\mu\big),\qquad
\forall u\in C^\infty_0(K).
\end{equation}
Our statement is now a consequence of a classical localization
procedure, that we briefly recall here for the sake of
completeness. Without loss of generality, we can assume that $A$
is properly supported and, by homogeneity, that $|\eta_0|=1$. Take
$\chi \in  C^\infty_0(\R^\nu)$ with $\chi\equiv 1$ near $K \cup\
\textrm{supp}\,Au$ for any $u\in C^\infty_0(K)$; since $Au=\chi
A(\chi u)$ we can assume that $a(y,\eta)$ is compactly supported
in $y$. Fix a $v\in C^\infty_0(\R^\nu)$ and put (with $t\geq 1$)
\begin{equation}\label{funzloc}
u_t(y)=e^{it^{2}\<y,\eta_0 \>}v(t(y-y_0))
\end{equation}
For $t$ large, $u_t \in C^\infty_0(K)$ and, after a few
computations, one gets:
\begin{equation}\label{eqn17}
\ds \left\{ \begin{array}{l}
         \widehat{u}_t(\eta)=\ds t^{-\nu}e^{i\<y_0,t^{2}\eta_0-\eta\>}
         \widehat{v}\big(\frac{\eta}{t}-t\eta_0\big) \\
         \vspace{-0.3cm}\\
         \parallel u_t\parallel^2_s=
         t^{4s-\nu}
         \big( \parallel v \parallel^2_{0}+ o(1) \big)\;,\;\;\;s\in\R
         \end{array}
\right.
\end{equation}
On the other hand,
$$Au_t(y)=e^{it^{2}\<y,\eta_0\>}\phi_t\big( t(y-y_0)\big)$$
with
$$\ds \phi_t(y)=(2\pi)^{-\nu} \int e^{i\<y,\eta\>}a\Big(y_0+\frac{y}{t},
t\eta+t^{2}\eta_0\Big)\widehat{v}(\eta) d\eta\;.$$ An application
of the Taylor formula yields
$$
a\big(y_0+\frac{y}{t},t\eta+t^{2}\eta_0\big)=
 a_{m'}\big(y_0+\frac{y}{t},\frac{\eta}{t}+\eta_0\big)
t^{2m'}+a_{m'-1/(h+1)}\big(y_0+\frac{y}{t},\frac{\eta}{t}+\eta_0\big)
t^{2m'-\frac{2}{h+1}}+O\big(t^{2m'-\frac{4}{h+1}}\big)
$$
\begin{equation*}
\hspace{-2,1cm} =a_{m'-1/(h+1)}\big(y_0,\eta_0\big)
t^{2m'-\frac{2}{h+1}}+O(t^{2m'-\g}),
\end{equation*}
where $\g=\min\{1,4/(h+1)\}$. Thus, by using inequality
(\ref{pal}) with $s=0$ and the estimate in (\ref{eqn17}), we obtain
$$\ds  \parallel Au_t \parallel^2_{0} = t^{4m'-\frac{4}{h+1}-\nu}|a_{m'-1/(h+1)}(y_0,\eta_0)|^2\|v\|_0^2+O(t^{4m'-\frac{2}{h+1}-\g-\nu})$$
\begin{equation}\label{disuguaI}
\qquad \geq \frac{1}{C} t^{4m'-4\sigma-\nu}\big(
\parallel v \parallel^2_{0}+o(1) \big)
-t^{4\mu-\nu}\big(\parallel v \parallel^2_{0}+o(1) \big)
\end{equation}
as $t\rightarrow +\infty$. Since $h>1$ one has $\g>2/(h+1)$, so
that
$O(t^{4m'-\frac{2}{h+1}-\g-\nu})=o(t^{4m'-\frac{4}{h+1}-\nu})$.
Thus dividing by $t^{4m'-\frac{4}{h+1}-\nu}$, and letting
$t\rightarrow +\infty$ shows that the left-hand-side is bounded.
As a consequence, we must have
$\sigma\geq\frac{1}{h+1}.$\\
Moreover, it is worth noting that the inequality \eqref{pal} holds for $\sigma=\frac{1}{h+1}$ only if $a_{m'-1/(h+1)}(y_0,\eta_0)\neq 0$.\\
Assume now that $a_{m'-\frac{1}{h+1}}(y_0,\eta_0)=0$ and $h$ is
odd (hence $h\geq 3$). In this case we get
$$
a\big(y_0+\frac{y}{t},t\eta+t^{2}\eta_0\big)=a_{m'-2/(h+1)}\big(y_0,\eta_0\big)
t^{2m'-\frac{4}{h+1}}+O(t^{2m'-\g}),
$$
with $\g=\min\{1,6/(h+1)\}\geq 4/(h+1)$ (note that $\g> 4/(h+1)$
if $h>3$). Inserting the above relation in (\ref{pal}) we have
$$\ds  \parallel Au_t \parallel^2_{0} = t^{4m'-\frac{8}{h+1}-\nu}|a_{m'-2/(h+1)}(y_0,\eta_0)|^2\|v\|_0^2+O(t^{4m'-\frac{4}{h+1}-\g-\nu})$$
\begin{equation}\label{disuguaII}
\qquad \geq \frac{1}{C} t^{4m'-4\sigma-\nu}\big(
\parallel v \parallel^2_{0}+o(1) \big)
-t^{4\mu-\nu}\big(\parallel v \parallel^2_{0}+o(1) \big)
\end{equation}
and letting again $t\rightarrow +\infty$ yields
$\sigma\geq\frac{2}{h+1}.$
We complete the proof by observing that, if $h>3$, the inequality \eqref{pal} holds for $\sigma=\frac{2}{h+1}$ only if $a_{m'-2/(h+1)}(y_0,\eta_0)\neq 0$. \fine \\
As a by-product of the proof above we get the following remark.
\begin{remark}\label{ossnec}
$A$ is hypoelliptic  with the minimal loss  of derivatives
$\sigma=\frac{1}{h+1}$ only if\linebreak
$a_{m'-1/(h+1)}(y,\eta)\neq 0$ whenever $a_{m'}(y,\eta)=0$. On the
other hand, if $h\neq3$ and \linebreak
$a_{m'-1/(h+1)}(y_0,\eta_0)=0$, the minimal loss
$\sigma=\frac{2}{h+1}$ can be attained only if
$a_{m'-2/(h+1)}(y_0,\eta_0)\neq 0$. The case $h=3$  will be
discussed later on. As a matter of fact, the case $h=3$ deserves a
particular analysis: in this setting, via Lemma \ref{ine} the
hypoellipticity of $A$ means that $A$ satisfies the following
inequality (for any $\mu<m'-\frac{1}{2}$)
$$\langle A^\ast Au, u\rangle=\|Au\|_0^2 \geq c\|u\|^2_{m'-\frac{1}{2}}-C\|u\|_\mu^2,\qquad \forall u\in C^\infty_0(K).$$
Thus all the terms of order greater o equal to $2m'-1$ in the
symbol of $\sigma(A^\ast A)$ matter, in particular the
$(y_j,\eta_j)$-derivatives of the principal symbol
$a_{m'}(y,\eta)$ of $A$. This happens only in the case $h=3$, in
the other cases the above derivatives are negligible.
\end{remark}
Our aim is now to find out necessary and sufficient conditions
ensuring the hypoellipticity of $A$ with the minimal loss of
derivatives shown in Lemma \ref{loss}. Since this result will be
applied to the operator $L$ in Theorem \ref{main}, from now on we
assume that $a_{m'-1/(h+1)}\equiv 0$ whenever $h$ is an odd
integer (see Remark \ref{struttura}). In view of Remark
\ref{ossnec}, we firstly assume $h\neq 3$.
\begin{proposition}\label{cns}
Let $h$ be a positive integer with $h\neq3$. Assume that for every
$(y_0,\eta_0)\in T^\ast\mathbb{S}^\nu: a_{m'}(y_0,\eta_0)=0$,
there exist a  neighborhood $U$ of $(y_0,\eta_0)$ in
$T^\ast\mathbb{S}^\nu$ and a  constant $c=c(U)>0$ such that:
\begin{itemize}
\item[$(H1).$] $a_{m'-\frac{r}{h+1}}(y_0,\eta_0)\neq 0$;
\item[$(H2).$] for every $(y,\eta)\in U$ one has $$-{\rm{Re}}\
a_{m'}a_{m'-\frac{r}{h+1}}(y,\eta)\leq
|a_{m'}(y,\eta)|\sqrt{|a_{m'-\frac{r}{h+1}}(y,\eta)|^2-c},$$
\end{itemize}
where $r=1$ if $h$ is even, whereas $r=2$ if $h$ is odd.\\
Then $A$ is $C^\infty$ hypoelliptic with loss of $r/(h+1)<1$ derivatives.\\
On the other hand, if $A$ is $C^\infty$ hypoelliptic with loss of
$r/(h+1)$ derivatives, then $A$ verifies the hypotheses $(H1)$ and
$(H2)$.
\end{proposition}
% la condizione H2 si puo' realizzare anche se -{\rm{Re}}\ a_{m'}a_{m'-\frac{r}{h+1}}>0. Infatti tale condizione si riduce a dire che
% Re A*B \leq |A|*|B| -c|A|^2. Scegliendo ad esempio (utile la rappresentazione goniometrica per capire il fenomeno)
% A=\sqrt{2}+i\sqrt{2} B=1+i\sqrt{3} ------> \sqrt{2}-\sqrt{6}<4-4c con 0<c<<1
{\bf Proof.} We start off by showing that $(H1)$ and $(H2)$ are sufficient conditions for the hypoellipticity of $A$. To this end we shall construct a parametrix of
$A$ by using Theorem 22.1.3 Vol.III \cite{ho}.
%In view of Lemma \ref{ine} it is enough to prove inequality \eqref{stimaauto} with $\sigma=r/(h+1)<1$
Consider, for any bounded open set $\Omega\subset\R^\nu$ and every
$(y,\eta)\in T^\ast\Omega\setminus 0, |\eta|\geq1$,
$$|a_{m'}(y,\eta)+a_{m'-\frac{r}{h+1}}(y,\eta)|^2=|a_{m'}(y,\eta)|^2+2{\rm{Re}}\  a_{m'}a_{m'-\frac{r}{h+1}}(y,\eta)+|a_{m'-\frac{r}{h+1}}(y,\eta)|^2$$
$$=|\eta|^{2m'-2\frac{r}{h+1}}\big(\underbrace{|a_{m'}(y,\omega)|^2t^2+2{\rm{Re}}\  a_{m'}a_{m'-\frac{r}{h+1}}(y,\omega)t+|a_{m'-\frac{r}{h+1}}(y,\omega)|^2}_{=g(y,\omega)(t)}\big),$$
with $\omega=\eta/|\eta|$ and $t=|\eta|^{\frac{r}{h+1}}$. As a
consequence of hypotheses $(H1)$ and $(H2)$ the non-negative
function $g(y,\omega)(t)$ is actually bounded from below by a
positive constant depending on $\Omega$. Precisely, if
$a_{m'}(y,\eta)=0$ this is due to $(H1)$; in the other cases,
$g(y,\omega)(t)$ is a parabola in the $t$-variable. Thus
conditions $(H1)$ and $(H2)$ ensure that, conically near the
characteristic set $\textrm{Char}(A)$ of $A$, its minimum in the
region $t\geq0$ is uniformly bounded from below by a positive
constant; whereas, conically outside $\textrm{Char}(A)$, this is a
trivial consequence of the ellipticity of $A$. Therefore, for
every $(y,\eta)\in T^\ast\Omega\setminus 0, |\eta|\geq1$, we have
$$|a_{m'}(y,\eta)+a_{m'-\frac{r}{h+1}}(y,\eta)|\geq c'|\eta|^{m'-\frac{r}{h+1}},$$
whence,  for every $(y,\eta)\in T^\ast\Omega\setminus 0,
|\eta|>>1$,
$$|a(y,\eta)|\geq |a_{m'}(y,\eta)+a_{m'-\frac{r}{h+1}}(y,\eta)|-|\big(a-\sum_{j=0}^r a_{m'-\frac{j}{h+1}}\big)(y,\eta)| \geq c''|\eta|^{m'-\frac{r}{h+1}}.$$
Since $a(y,\eta)\in S_{1,0}^{m'}$, we obtain
$$|\partial_y^\a\partial_\eta^\b a(y,\eta)|\leq C_{\a,\b}(1+|\eta|)^{m'-\frac{r}{h+1}}(1+|\eta|)^{\frac{r}{h+1}-|\beta|}\leq C'_{\a,\b}|a(y,\eta)|
(1+|\eta|)^{\delta|\alpha|-\rho|\beta|},$$
where $\rho=1-\frac{r}{h+1}>\delta=\frac{r}{h+1}$ if $h\neq3$.\\
Finally from Lemma 22.1.2 and Theorem 22.1.3 Vol.III \cite{ho} we
can construct
a parametrix of $A$ in the H\"ormander class ${\rm{OP}}S_{\rho,\delta}^{-m'+\frac{r}{h+1}}$, which proves the first part of Theorem \ref{cns}.\\
The necessity of the condition $(H1)$ follows immediately from Remark \ref{ossnec}. As for condition $(H2)$ we proceed by localizing the estimate \eqref{stimaauto} with $\sigma=\frac{r}{h+1}$ and $s=0, \mu<m'-\frac{r}{h+1}$. Here we assume that $h$ is even; we shall see that a similar argument applies to the odd case.\\
From  Lemma \ref{ine} we get
\begin{equation}\label{non}
\langle A^\ast Au, u\rangle=\|Au\|_0^2 \geq
c\|u\|^2_{m'-\frac{1}{h+1}}-C\|u\|_\mu^2,\qquad \forall u\in
C^\infty_0(K).
\end{equation}
This estimate is stable under perturbation of order less than
$\nu<m'-\frac{1}{h+1}$, therefore we can assume that
$a_{m'-\frac{j}{h+1}}\equiv 0$ if $j>1$. A straightforward
computation yields
\begin{equation}\label{forma}
\sigma(A^\ast A)=|a_{m'}|^2+2{\rm Re}\
(a_{m'}a_{m'-\frac{1}{h+1}})+|a_{m'-\frac{1}{h+1}}|^2\qquad {\rm
mod.}\ \, S^{2m'-1}.
\end{equation}
Let $(\overline{y},\overline{\eta})\in T^\ast\R^\nu\setminus 0$
with $|\overline{\eta}|=1$ and consider the localizing function
\eqref{funzloc} with $(y_0,\eta_0)$ replaced by
$(\overline{y},\overline{\eta})$
\begin{equation}\label{localfunction}
u_t(x)=e^{it^{2}\<y,\overline{y} \>}v(t(y-\overline{y})).
\end{equation}
Here we assume that $0\neq v\in C^\infty_0(\R^\nu)$ is an even
function. By arguing as in the proof of Lemma \ref{loss} we obtain
$$\langle A^\ast Au_t, u_t\rangle=t^{-\nu}\sigma(A^\ast A)(\overline{y},t^2\overline{\eta})\|v\|_0^2+t^{-\nu-1}\sum_{j=1}^\nu\partial_{y_j}\sigma(A^\ast A)(\overline{y},t^2\overline{\eta})\cdot\langle y_j v, v\rangle$$
$$+t^{-\nu+1}\sum_{j=1}^\nu\partial_{\xi_j}\sigma(A^\ast A)(\overline{y},t^2\overline{\eta})\cdot\langle D_j v, v\rangle+O(t^{4m'-\nu-2}).$$
Due to the compact support of $v(x)$ we have $\langle D_j v,
v\rangle=0$ and, furthermore, from its parity we get $\langle y_j
v, v\rangle=0$; thus, by \eqref{forma} we obtain
$$\langle A^\ast Au_t, u_t\rangle=t^{-\nu}\sigma(A^\ast A)(\overline{y},t^2\overline{\eta})\|v\|_0^2+O(t^{4m'-\nu-2}).$$
whence,
$$\langle A^\ast Au_t, u_t\rangle=t^{-\nu}\|v\|_0^2\Big(|a_{m'}(\overline{y},\overline{\eta})|^2t^{4m'}+
2{\rm Re}\
(a_{m'}a_{m'-\frac{1}{h+1}})(\overline{y},\overline{\eta})t^{4m'-\frac{2}{h+1}}$$
$$+|a_{m'-\frac{1}{h+1}}(\overline{y},\overline{\eta})|^2 t^{4m'-\frac{4}{h+1}}\Big)
+o(t^{4m'-\frac{4}{h+1}-\nu}),$$ due to the fact that
$4m'-\nu-2<4m'-\nu-\frac{4}{h+1}$. By inserting the above relation
into \eqref{non} and using \eqref{eqn17} we get
$$\Big(|a_{m'}(\overline{y},\overline{\eta})|^2t^{4m'}+
2{\rm Re}\
(a_{m'}a_{m'-\frac{1}{h+1}})(\overline{y},\overline{\eta})t^{4m'-\frac{2}{h+1}}
+|a_{m'-\frac{1}{h+1}}(\overline{y},\overline{\eta})|^2 t^{4m'-\frac{4}{h+1}}\Big)\|v\|_0^2$$
\begin{equation}\label{dico}
\geq\frac{1}{C} t^{4m'-\frac{4}{h+1}}
\parallel v \parallel^2_{0}+
o(t^{4m'-\frac{4}{h+1}}).
\end{equation}
Notice that here the error term $o(t^{4m'-\frac{4}{h+1}})$ and the positive constant $1/C$ are uniform if $(\overline{y},\overline{\eta})$ takes values in a small neighborhood $U_0$ of $(y_0,\eta_0)$ in $T^\ast\mathbb{S}^\nu$.\\
Dividing  the above expression by $t^{4m'-\frac{4}{h+1}}\|v\|_0^2$
we get, for a suitable small positive constant $c_1$,
\begin{equation}\label{crucform}
|a_{m'}(y,\eta)|^2 (t^{\frac{2}{h+1}})^2+ 2{\rm Re}\
(a_{m'}a_{m'-\frac{1}{h+1}})(y,\eta)t^{\frac{2}{h+1}}
+|a_{m'-\frac{1}{h+1}}(y,\eta)|^2 \geq c_1>0,
\end{equation}
whenever  $(y,\eta)\in U_0$ and $t\geq c_0$ for a suitable large positive constant $c_0=c_0(U_0)$.\\
This fact suggests to consider the following parabola in the
variable $z\in\R$
$$q(z):=|a_{m'}(y,\eta)|^2 z^2+
2{\rm Re}\ (a_{m'}a_{m'-\frac{1}{h+1}})(y,\eta)z
+|a_{m'-\frac{1}{h+1}}(y,\eta)|^2.$$ It is enough to show that
$(H2)$ holds if ${\rm{Re}}\
a_{m'}a_{m'-\frac{1}{h+1}}(y,\eta)<0$, for otherwise
$(H2)$ is automatically satisfied.\\
In this case, $p(z)$ takes its minimum in $\ds
z_{\min}(y,\eta):=-\frac{{\rm{Re}}\
a_{m'}a_{m'-\frac{1}{h+1}}}{|a_{m'}|^2}(y,\eta)>0$ and
$$q\big(z_{\min}(y,\eta)\big)=-\frac{\Delta}{4|a_m|^2}(y,\eta)\geq 0,$$
with $\Delta(y,\eta)=4({\rm{Re}}\  a_{m'}a_{m'-\frac{1}{h+1}}(y,\eta))^2-4|a_{m'}(y,\eta)a_{m'-\frac{1}{h+1}}(y,\eta)|^2\geq 0$.\\
We observe that condition $(H2)$ amounts to saying that $p\big(z_{\min}(y,\eta)\big)\geq c$ for every $(y,\eta)$ in a neighborhood $U$ of $(y_0,\eta_0)$ in $T^\ast\mathbb{S}^\nu$ and we proceed by contradiction of the last statement.\\
If this were false then we could select a sequence
$(x_j,\omega_j)\in T^\ast\mathbb{S}^\nu$ such that
\begin{itemize}
\item[a)] $(x_j,\omega_j){}_{\overrightarrow{j\rightarrow
+\infty}} (y_0,\eta_0);$ \item[b)] ${\rm Re}\,
a_{m'}a_{m'-\frac{1}{h+1}}(x_j,\omega_j)<0$ for every $j\in\N$;
\item[c)] $\ds
q\big(z_{\min}(x_j,\omega_j)\big)=-\frac{\Delta}{4|a_m'|^2}(x_j,\omega_j)\;{}_{\overrightarrow{j\rightarrow
+\infty}}\;0.$
\end{itemize}
From $c)$ it follows that $\frac{\big({\rm Re}\,
a_{m'}a_{m'-\frac{1}{h+1}}\big)^2(x_j,\omega_j)}{|a_{m'}(x_j,\omega_j)|^2}
\;{}_{\overrightarrow{j\rightarrow
+\infty}}\;|a_{m'-\frac{1}{h+1}}(y_0,\eta_0)|^2>0$, so that
$$\left\{\begin{array}{l}
\ds z_{\min}(x_j,\omega_j)=-\frac{1}{{\rm{Re}}\  a_{m'}a_{m'-\frac{1}{h+1}}(x_j,\omega_j)}\cdot\frac{\big({\rm Re}\, a_{m'}a_{m'-\frac{1}{h+1}}\big)^2}{|a_{m'}|^2}(x_j,\omega_j)\;{}_{\overrightarrow{j\rightarrow +\infty}}\;+\infty,\\
q\big(z_{\min}(x_j,\omega_j)\big)\;{}_{\overrightarrow{j\rightarrow
+\infty}}\;0.
\end{array}\right.,$$
but this clearly contradicts \eqref{crucform}, thus the proof of condition $(H2)$ is complete if $h$ is even.\\
Let us briefly discuss the odd case $h>3$. We are going to show
that if $A$ is $C^\infty$ hypoelliptic with loss of $2/(h+1)$
derivatives, then the hypothesis $(H2)$ holds true. Again via
Lemma \ref{ine} the hypoellipticity of $A$ is equivalent to the
following inequality (for any $\mu<m'-\frac{2}{h+1}$)
\begin{equation}\label{stimadalbassoII}
\langle A^\ast Au, u\rangle=\|Au\|_0^2 \geq
c\|u\|^2_{m'-\frac{2}{h+1}}-C\|u\|_\mu^2,\qquad \forall u\in
C^\infty_0(K).
\end{equation}
Notice that the above estimate is stable if we add to $A$ any pseudodifferential operator of order less than $m'-\frac{2}{h+1}$. Therefore, without loss of generality, we can assume that the terms $a_{m'-\frac{r}{h+1}}=0,\; r>2,$ in the asymptotic expansion of the symbol of $A$ identically vanish.\\
By using the localizing function \eqref{localfunction} and
proceeding as above we get
$$\langle A^\ast Au_t, u_t\rangle=t^{-\nu}\sigma(A^\ast A)(\overline{y},t^2\overline{\eta})\|v\|_0^2+O(t^{4m'-\nu-2}).$$
A direct computation shows that
\begin{equation}\label{forma2}
\sigma(A^\ast A)=|a_{m'}|^2+2{\rm Re}\
(a_{m'}a_{m'-\frac{2}{h+1}})+|a_{m'-\frac{2}{h+1}}|^2\qquad {\rm
mod.}\ \, S^{2m'-1},
\end{equation}
whence it readily follows that
$$\langle A^\ast Au_t, u_t\rangle=t^{-\nu}\|v\|_0^2\Big(|a_{m'}(\overline{y},\overline{\eta})|^2t^{4m'}+
2{\rm Re}\
(a_{m'}a_{m'-\frac{2}{h+1}})(\overline{y},\overline{\eta})t^{4m'-\frac{4}{h+1}}$$
$$+|a_{m'-\frac{2}{h+1}}(\overline{y},\overline{\eta})|^2 t^{4m'-\frac{8}{h+1}}\Big)
+o(t^{4m'-\frac{8}{h+1}-\nu}),$$ due to the fact that
$4m'-\nu-2<4m'-\nu-\frac{8}{h+1}$ if $h>3$. Inserting this formula
in \eqref{stimadalbassoII} yields
$$\Big(|a_{m'}(\overline{y},\overline{\eta})|^2 t^{4m'}+
2{\rm Re}\
(a_{m'}a_{m'-\frac{2}{h+1}})(\overline{y},\overline{\eta})t^{4m'-\frac{4}{h+1}}
+|a_{m'-\frac{2}{h+1}}(\overline{y},\overline{\eta})|^2 t^{4m'-\frac{8}{h+1}}\Big)\|v\|_0^2$$
\begin{equation}\label{dico2}
\geq\frac{1}{C} t^{4m'-\frac{8}{h+1}}
\parallel v \parallel^2_{0}+
o(t^{4m'-\frac{8}{h+1}}).
\end{equation}
and finally we have, for a suitable positive constant $c_1$,
$$
|a_{m'}(y,\eta)|^2(t^{\frac{4}{h+1}})^2+ 2{\rm Re}\
(a_{m'}a_{m'-\frac{1}{h+1}})(y,\eta)t^{\frac{4}{h+1}}
+|a_{m'-\frac{1}{h+1}}(y,\eta)|^2\geq c_1>0,
$$
whenever  $(y,\eta)\in U_0$ and $t\geq c_0$ for a suitable big
positive constant $c_0=c_0(U_0)$. From now on the proof proceeds
exactly as in the even case.\fine

Let us now consider the case $h=3$. In view of the examples we are interested in, we assume that $m'\leq 1$ (as a matter of fact, in Thm. \ref{teoprinc}  we need $m'=1/2$) and we limit ourselves to find out only sufficient conditions ensuring the
minimal loss of derivatives for $A$.
\begin{proposition}\label{h=3}
Let $h=3$. Assume that for every $(y_0,\eta_0)\in
T^\ast\mathbb{S}^\nu: a_{m'}(y_0,\eta_0)=0$, there exist a
neighborhood $U$ of $(y_0,\eta_0)$ in $T^\ast\mathbb{S}^\nu$ and a
constant $c=c(U)>0$ such that:
\begin{itemize}
\item[$(H1).$] $|a_{m'-\frac{1}{2}}(y_0,\eta_0)|^2+\frac{1}{i}\{ \overline{a_{m'}},a_{m'}\}(y_0,\eta_0)>0$; \item[$(H2).$] for every $(y,\eta)\in U$ one has
$$-{\rm{Re}}\ a_{m'}a_{m'-\frac{1}{2}}(y,\eta)\leq |a_{m'}(y,\eta)|\sqrt{|a_{m'-\frac{1}{2}}(y,\eta)|^2+\frac{1}{i}\{ \overline{a_{m'}},a_{m'}\}(y,\eta)-c},$$
where $\{f,g\}$ denotes the Poisson brackets of the functions $f,g$.
\end{itemize}
Then $A$ is $C^\infty$ hypoelliptic with the minimal loss of $2/(h+1)=1/2$ derivatives.\\
%Furthermore, if the principal symbol $a_{m'}$ vanishes identically on $U$, then $(H1)$ is actually a necessary and sufficient condition to have such a loss of regularity.
Furthermore, if the principal symbol $a_{m'}(y,\eta)$ vanishes to the second order at $(y,\eta)=(y_0,\eta_0)$, then $(H1)$ reads as $a_{m'-\frac{1}{2}}(y_0,\eta_0)\neq 0$ and is actually a necessary and sufficient condition to have such a loss of regularity.
\end{proposition}
{\bf Proof.} We observe that $A=A'+A''$ with
$A'=\textrm{Op}(a_{m'}+a_{m'-\frac{1}{2}})$ and
$A''\in\textrm{OP}S^{m'- \frac{3}{4}}$ and thus, from the Sobolev
continuity of $A''$, we have
$$\|Au\|_0 \geq \|A'u\|_0-\|A''u\|_0\geq \|A'u\|_0 - C\|u\|_{m'- \frac{3}{4}},\quad \forall u\in C^\infty_0(K).$$
Therefore, in view of Lemma \ref{ine}, it suffices to
show that
$$\langle A'^*A'u,u\rangle=\|A'u\|_0^2 \geq c\|u\|_{m'-\frac{1}{2}}^2-C\|u\|_{m'-1}^2,\qquad \forall u\in C^\infty_0(K),$$
 $c$ being a positive constant.\\
This task will be achieved by applying the Fefferman-Phong
inequality to the operator $A'^\ast A'$; to this aim it is
convenient to use the Weyl quantization (see Section 18.5 Vol.III
\cite{ho}). By Theorems 18.5.4 and 18.5.10 Vol.III \cite{ho} we
have $\sigma^w(A'^*)=\overline{\sigma^w(A')}$ and
$$\sigma^w(A')=a_{m'}(y,\eta)+a_{m'-\frac{1}{2}}(y,\eta)+i/2\sum_{j=1}^\nu
\partial_{\eta_j}\partial_{y_j} a_{m'}(y,\eta)\quad \textrm{mod.}\,S^{m'- \frac{3}{2}},$$
whence
\begin{eqnarray*}
\sigma^w(A'^\ast
A')(y,\eta)&=&|a_{m'}(y,\eta)+a_{m'-\frac{1}{2}}(y,\eta)|^2+
\frac{1}{i}\{\overline{a_{m'}},a_{m'}\}(y,\eta)\\
& & +\textrm{Im}\;a_{m'}(y,\eta)\sum_{j=1}^\nu
\partial_{\eta_j}\partial_{y_j} a_{m'}(y,\eta)+e(y,\eta)\\
&=&|a_{m'}(y,\eta)|^2+2\textrm{Re}\;a_{m'}(y,\eta)a_{m'-\frac{1}{2}}(y,\eta)+|a_{m'-\frac{1}{2}}(y,\eta)|^2\\
& &
+\frac{1}{i}\{\overline{a_{m'}},a_{m'}\}(y,\eta)+\textrm{Im}\;a_{m'}(y,\eta)\sum_{j=1}^\nu
\partial_{\eta_j}\partial_{y_j} a_{m'}(y,\eta)+e(y,\eta),
\end{eqnarray*}
where $e(y,\eta)$ is an error term of order $2m'- \frac{3}{2}$, i.e. $e(y,\eta)\in S^{2m'- \frac{3}{2}}(\R^\nu\times\R^\nu)$. Notice that $\sigma^w(A'^\ast A')(y,\eta)$ is a real symbol,  $A'^\ast A'$ being a formally selfadjoint operator.\\
By arguing as in the proof of Proposition \ref{cns}, from the
hypotheses $(H1)$, $(H2)$ of Proposition \ref{h=3} it follows
that, conically near the characteristic set $\Sigma=\{(y,\eta) |
a_{m'}(y,\eta)=0,\, \eta\neq 0\}$,
$$|a_{m'}(y,\eta)|^2+2\textrm{Re}\;a_{m'}(y,\eta)a_{m'-\frac{1}{2}}(y,\eta)+|a_{m'-\frac{1}{2}}(y,\eta)|^2$$
$$
+\frac{1}{i}\{\overline{a_{m'}},a_{m'}\}(y,\eta)\geq
\frac{1}{c_1}|\eta|^{2m'-1},\quad y\in K\subset\subset \R^\nu,\;\;
|\eta|\geq c_1,$$ for a positive large constant $c_1=c_1(K)$.
Since $e(y,\eta)\in S^{2m'- \frac{3}{2}}$ and the term $S^{2m'-
1}\ni\textrm{Im}\;a_{m'}\sum_{j=1}^\nu
\partial_{\eta_j}\partial_{y_j} a_{m'}$ vanishes on $\Sigma$, we easily see that, conically near $\Sigma$,
$$\sigma^w(A'^\ast A')(y,\eta)\geq \frac{1}{c_2}|\eta|^{2m'-1},\quad y\in K\subset\subset \R^\nu,\;\; |\eta|\geq c_2,$$
for a new constant $c_2=c_2(K)>0$. In view of the
ellipticity of $A'$ outside $\Sigma$, we finally have that
$$\sigma^w(A'^\ast A')(y,\eta)-\frac{1}{c_3}|\eta|^{2m'-1}\geq 0,\quad y\in K\subset\subset \R^\nu,\;\; |\eta|\geq c_3,\quad c_3=c_3(K)>0.$$
An application of the Fefferman-Phong inequality (see
also Theorem 18.6.8 and Corollary 18.6.11 Vol.III \cite{ho}) to
the operator with symbol $\sigma^w(A'^\ast
A')(y,\eta)-\frac{1}{c_3}|\eta|^{2m'-1}\in S^{2m'}$ gives
$$\|A'u\|_0^2=\langle A'^*A'u,u\rangle \geq \frac{1}{c_3}\|u\|_{m'-\frac{1}{2}}^2-C\|u\|_{m'-1}^2,\qquad \forall u\in C^\infty_0(K).$$
This completes the first part of the statement of the proposition \ref{h=3}.\\
Assume now that the principal symbol $a_{m'}$ to the second order at $(y_0,\eta_0)$; we are left to show that $(H1)$ is also a necessary condition for the minimal hypoellipticity of $A$.\\
We argue similarly as done in the second part of the proof of
Prop. \ref{cns}. By using the localizing function $u_t(x)$ in
\eqref{localfunction} and Lemma \ref{ine} we have, for any
$\mu<m'-\frac{1}{2}$,
$$\langle A^\ast Au_t, u_t\rangle=\|Au_t\|_0^2 \geq c\|u_t\|^2_{m'-\frac{1}{2}}-C\|u_t\|_\mu^2.$$
Since $a_{m'}(y,\eta)$ vanishes to the second order at $(y_0,\eta_0)$, we get
\begin{eqnarray*}
\langle A^\ast Au_t, u_t\rangle&=&t^{-\nu}\sigma(A^\ast A)(y_0,t^2\eta_0)\|v\|_0^2+o(t^{4m'-\nu-2}),\\
&=&t^{-\nu}\|v\|_0^2|a_{m'-\frac{1}{2}}(y_0,\eta_0)|^2t^{4m'-2}
+o(t^{4m'-2-\nu}).
\end{eqnarray*}
Due to \eqref{eqn17} and by inserting this formula in the above
inequality, we obtain
$$|a_{m'-\frac{1}{2}}(y_0,\eta_0)|^2t^{4m'-2}\|v\|_0^2
+o(t^{4m'-2})\geq C t^{4m'-2}\parallel v
\parallel^2_{0}+o(t^{4m'-2}).$$ By letting $t\longrightarrow
+\infty$, we complete the proof.\fine\\
We end this section by analyzing the hypoellipticity of $A=a(y,D)$ in \eqref{amj} a special setting.
\begin{proposition}\label{tangential}
Assume that $a_{m'}(y_0,\eta_0)=0$ for some $(y_0,\eta_0)\in T^\ast\mathbb{S}^\nu$ and suppose that $a_{m'-\frac{j}{h+1}}\equiv 0$ for every $0<j<h+1$. Then $A$ can
be hypoelliptic only with loss of $\delta\geq 1/2$ derivatives. Furthermore, the threshold is realized if and only if
\begin{equation}\label{popa}
\frac{1}{i}\{a_{m'},\overline{a_{m'}}\}(y_0,\eta_0)<0.
\end{equation}
\end{proposition}
{\bf Proof.} In order to prove the first part of the statement we could argue as done in Lemma \ref{loss}. However, due to the fact that
$A-a_{m'}(y,D_y)\in \textrm{OP}S^{m'-1}$, from Lemma \ref{ine} it is easily seen that the hypoellipticity with loss of $\delta<1$ derivatives only relies on the principal symbol $a_{m'}(y,\eta)$. Hence, by applying Prop. 27.1.8 in H\"{o}rmander Vol.IV \cite{ho} (setting $n=\nu$, $k=1$) one readily gets $\delta\geq k/(k+1)=1/2$. Furthermore, from \eqref{popa} it follows that $H_{\textrm{Re}\, a_{m'}}\textrm{Im}\, a_{m'}(y_0,\eta_0)=\{\textrm{Re}\, a_{m'}, \textrm{Im}\, a_{m'}\}(y_0,\eta_0)=-\frac{1}{2i}\{a_{m'},\overline{a_{m'}}\}(y_0,\eta_0)>0$, whence an application of Th. 27.1.11 in H\"{o}rmander Vol.IV \cite{ho} completes the proof.

\renewcommand\thesection{\Alph{section}}
\section{Appendix II. The $h$-pseudodifferential calculus}\label{calc}

\setcounter{equation}{0}
\setcounter{theorem}{0}

\renewcommand{\thetheorem}{\thesection.\arabic{theorem}}
\renewcommand{\theproposition}{\thesection.\arabic{proposition}}
\renewcommand{\thelemma}{\thesection.\arabic{lemma}}
\renewcommand{\thedefinition}{\thesection.\arabic{definition}}
\renewcommand{\thecorollary}{\thesection.\arabic{corollary}}
\renewcommand{\theequation}{\thesection.\arabic{equation}}

This appendix is devoted to the construction of the pseudodifferential calculus used  in the study of the hypoellipticity of the Grushin-type model \eqref{grushin}.
Following the ideas of \cite{boutet} we first define the classes of symbols we deal with, taking into account the anisotropy due to the parameter $h>1$ (see
\cite{bmt1}, \cite{bmt2} for further details).
\begin{definition}\label{definition11}
Let $m,k\in\R.$ By $S_h^{m,k}$ we denote the class of all smooth
functions $a(x,\xi) \colon \R^n\times\R^n\longrightarrow\C$ such
that, for all multi-indices $\alpha,\beta,\gamma,\delta$
\begin{equation}\label{eq21}
\left|\p_{x_1}^\alpha\p_{x'}^\beta\p_{\xi_1}^\gamma\p_{\xi'}^\delta
a(x,\xi)\right|\lesssim |\xi|^{m-\gamma-|\delta|}
\left(\frac{|\xi_1|}{|\xi|}+|x_1|^h+\frac{1}{|\xi|^{h/(h+1)}}\right)^{k-\alpha/h-\gamma}.
\end{equation}
We denote by $\OPS_h^{m,k}$ the corresponding class of
pseudodifferential operators. We set
$$S_h^{m,\infty}:=\bigcap_{k\in\R}S_h^{m,k}.$$
For $m,k\in\R$ we put
\begin{equation}
\mathcal{H}_h^{m,k}:=\bigcap_{j\geq 0}S_h^{m-j,k-j\frac{h+1}{h}},
\label{eq34}\end{equation} and $\mathrm{OP}\mathcal{H}_h^{m,k}$
denotes the corresponding class of operators.
\end{definition}
For instance, the operator \eqref{grushin2} belongs to $\OPS_h^{2,2}$.\\
By a straightforward computation, see e.g. \cite{marco}, we have that $S^{m,k}_h \subseteq{S^{m+\frac{h}{h+1}k_-}_{\frac{1}{h+1},\frac{1}{h+1}}}$ (with
$k_-=\max\{0,-k\}$). Moreover, an easy check shows that $m\leq{m'}$ and $\ds m-\frac{h}{h+1}k\leq{m'-\frac{h}{h+1}k'}$
then $S^{m,k}_h\subseteq{S^{m',k'}_h}$.\\
Let $ f_{-j}(x, \xi) \in S_{h}^{m, k+\frac{j}{h}} $, then there
exists $ f(x, \xi) \in S_{h}^{m, k} $ such that $ f \sim \sum_{j
\geq 0} f_{-j} $, i.e. $ f - \sum_{j=0}^{N-1} f_{-j} \in S_{h}^{m,
k + \frac{N}{h}} $, thus $ f $ is defined modulo a symbol
in $ S_{h}^{m, \infty}$.\\
The next class of symbols is introduced having in mind the structure of the localized operator \eqref{locgrushin}.
\begin{definition}\label{bbS}
Given $k \in \R$, we denote by $\mathbb{S}^{k}_h$ the space of the global symbols $b(x_1,\xi_1)\in C^\infty(\R_{x_1}\times \R_{\xi_1})$ such that, for  any mindex
$\alpha,\beta$, we have
\begin{equation}\label{migliore}
|\partial^\alpha_{x_1} \partial^\beta_{\xi_1}b(x',\xi')| \leq C (1+|x_1|^h+|\xi_1|) ^{k-\frac{\alpha}{h}-\beta}.
\end{equation}
Moreover, we define $\ds \mathbb{S}^{-\infty}_h=\bigcap_{k} \mathbb{S}^k_h$. We say that $b(x_1,\xi_1)$ is $h$-globally elliptic if $|b(x_1,\xi_1)|\geq
C(1+|x_1|^h+|\xi_1|)^k$. Finally we denote by $\mathrm{OP}\mathbb{S}^k_h$ (resp., in $\mathrm{OP}\mathbb{S}^{-\infty}_h$) the corresponding class of operators.
Throughout this paper we usually refer to them as the $h$-globally pseudodifferential operators.
\end{definition}
It is easily seen that  $P_\v$ (see \eqref{locgrushin}) is an operator in $\mathrm{OP}\S^2_h$, smoothly dependent on the parameter $\varrho=(x',\xi')$, and is $h$-globally elliptic.\\
Furthermore, we note that the  model operator \eqref{grushin} satisfies an intrinsic global homogeneous property that does not appear in the general class
$\OPS_h^{m,2}$ and is crucial in the calculus. Precisely, let use introduce the following definition.
\begin{definition}
\label{def:gh} We say that a symbol $ a(x, \xi) $ is globally
homogeneous (abbreviated g.h.) of degree $ m $ if, for $ \lambda
\geq 1 $, $ a(\lambda^{-1/(h+1)}x_1 , x', \lambda^{1/(h+1)}\xi_1,
\lambda \xi') = \lambda^{m} a(x,  \xi)$.
\end{definition}
Let $ f_{-j} $ be globally homogeneous of degree $ m - k
\frac{h}{h+1} - \frac{j}{h+1} $ and such that for every
multi-indices $ \gamma, \alpha, \beta$ satisfies the estimates
\begin{equation}
\label{defhom} \left|
\partial_{(x',\xi')}^{\gamma}\partial_{x_1}^{\alpha}
\partial_{\xi_1}^{\beta}
  f_{-j}(x, \xi) \right| \lesssim  \left(
  |\xi_1| + |x_1|^{h} + 1 \right)^{k- \frac{\alpha}{h}-\beta}, \qquad
(x_1, \xi_1) \in \R^{2},
\end{equation}
for $ (x', \xi') $ in a compact subset of $ \R^{n-1} \times
\R^{n-1}\setminus 0 $. Then $ f_{-j} \in S_{h}^{m, k +
\frac{j}{h}} $. For instance, the symbol $\sigma(P)$ of the
operator $P$ in \eqref{locgrushin} is g.h. of degree $2/(h+1)$ and
satisfies \eqref{defhom} for $k=2$, thus $\sigma(P)\in S_h^{2,2}$
as just observed.

In order to quantize the eigenfunctions $\phi_1, \phi_2$ in Lemma
\ref{quad}, we shall actually need the analog of the above
definition for functions not depending on $ \xi_1$.
\begin{definition}
\label{def:H} We denote by $ H_{h}^{m} $  the class of all smooth
functions  such that
\begin{equation}
\label{defsymb1} \left| \partial_{x'}^{\alpha}
\partial_{\xi'}^{\beta} \partial_{x_1}^{\gamma}
  a(x',\xi',x_1) \right| \lesssim |\xi'|^{m-j-|\beta|} \left(
  |x_1|^{h} + \frac{1}{|\xi'|^{\frac{h}{h+1}}} \right)^{-j\frac{h+1}{h}-
  \frac{\gamma}{h}}.
\end{equation}
Define the action of a symbol $ a(x',\xi',x_1) $ in $ H_{h}^{m} $
as the map $ a(x',D_{x'},x_1) \colon C_0^\infty(\R^{n-1}_{x'})
\longrightarrow C^\infty(\R^{n}_{x_1, x'}) $ defined by
$$
a(x',D_{x'},x_1) u (x_1, x') = (2\pi)^{-(n-1)} \int e^{ix'\xi'}
a(x',\xi',x_1) \hat{u}(\xi') d\xi'.
$$
Such an operator, modulo a regularizing operator (w.r.t. the $ t $
variable) is called a Hermite operator and we denote by $
OPH_{h}^{m} $ the
corresponding class.\\
We define the co-Hermite $ a^{*}(x',D_{x'},x_1) \colon
C_0^\infty(\R^n_{x_1, x'}) \longrightarrow C^\infty(\R^{n-1}_{x'})
$ as
$$
a^{*}(x',D_{x'},x_1) u (x') = (2\pi)^{-(n-1)} \int e^{ix'\xi'}
\overline{a(x',\xi', x_1)} \hat{u}(x_1, \xi') dx_1 d\xi'.
$$
We denote by $ {OPH_{h}^{*}}^{m} $ the related set of operators.
\end{definition}
Let $ \phi_{-j}(x',\xi',x_1) \in H_{h}^{m-\frac{j}{h+1}} $, then
there exists $ \phi(x',\xi',x_1) \in H_{h}^{m} $ such that $ \phi
\sim \sum_{j \geq 0} \phi_{-j} $, i.e. $ \phi - \sum_{j=0}^{N-1}
\phi_{-j} \in H_{h}^{m - \frac{N}{h+1}} $, so that $ \phi $ is
defined modulo a symbol
regularizing (w.r.t. the $ \xi' $ variables.)\\
Accordingly with Definition \ref{def:gh} we give the following
\begin{definition}
\label{def:gh2} We say that a symbol $ a(x', \xi', x_1) $ is
globally homogeneous (abbreviated g.h.) of degree $ m $ if, for
any $\lambda\geq 1$, one has that $
a(x',\lambda\xi',\lambda^{-1/(h+1)}x_1) =\lambda^{m} a(x', \xi',
x_1) $.
\end{definition}
The following remark is modelled on the functions $\phi_1, \phi_2$
in Lemma \ref{quad}.
\begin{remark}\label{ghherm}
Let $ \psi_{-j} $ be globally homogeneous of degree $ m -
\frac{j}{h+1} $ and such that for every multi-indices $\beta,
\alpha, \ell $ satisfies the estimates
\begin{equation}
\label{defhomherm} \left|
\partial_{(x',\xi')}^{\beta}\partial_{x_1}^{\alpha}
  \psi_{-j}(x, t,\tau) \right| \lesssim  \left(
  |x|^{h} + 1 \right)^{- \ell - \frac{\alpha}{h}}, \qquad x \in \R,
\end{equation}
for $ (x', \xi') $ in a compact subset of $ \R^{n-1} \times
\R^{n-1}\setminus 0 $. Then $ \psi_{-j} \in H_{h}^{m -
\frac{j}{h+1}} $.
\end{remark}
In the following lemma we compute the symbol of the formal adjoint
of operators in $OPH_{h}^{m}$, ${OPH_{h}^{*}}^{m}$ and in
$OPS_{h}^{m,k}$.
\begin{lemma}
\label{lemma:phi*} Let $ a \in H_{h}^{m} $, $ b \in S_{h}^{m,k} $;
then
\begin{itemize}
\item[(i)]{} the formal adjoint $ a(x', D_{x'}, x_1)^{*} $ belongs
to $ {OPH_{h}^{*}}^{m} $ and its symbol has the asymptotic
expansion
\begin{equation}
\label{eq:a*} \sigma(a(x', D_{x'}, x_1)^{*}) - \sum_{|\alpha|< N}
\frac{1}{\alpha!} \partial_{\xi'}^{\alpha} D_{x'}^{\alpha}
\overline{a(x', \xi', x_1)} \in H_{h}^{m - N}.
\end{equation}
\item[(ii)]{} The formal adjoint $ \left(a^{*}(x', D_{x'},
x_1)\right)^{*} $ belongs to $ OPH_{h}^{m} $ and its symbol has
the asymptotic expansion
\begin{equation}
\label{eq:a} \sigma(a^{*}(x', D_{x'}, x_1)^{*}) - \sum_{|\alpha|<
N} \frac{1}{\alpha!} \partial_{\tau}^{\alpha} D_{t}^{\alpha} a(x',
\xi', x_1) \in H_{h}^{m - N}.
\end{equation}
\item[(iii)]{} The formal adjoint $ b(x, D_{x})^{*} $ belongs to $
OPS_{h}^{m,k} $ and its symbol has the asymptotic expansion
\begin{equation}
\label{eq:b*} \sigma(a(x, D)^{*}) - \sum_{|\alpha|< N}
\frac{1}{\alpha!} \partial_{\xi}^{\alpha} D_{x}^{\alpha}
\overline{a(x, \xi)} \in S_{h}^{m - N, k - N \frac{h+1}{h}}.
\end{equation}
\end{itemize}
\end{lemma}
Next we give a brief description of the composition of the various
types of operator introduced so far.\\
As a matter of fact in the construction in Theorem \ref{main}  we
deal with asymptotic series of homogeneous symbols.
\begin{lemma}[\cite{h2}, Formula 2.4.9, \cite{marco} Prop. 1.11]
\label{lemma:SS} Let $ a \in S_{h}^{m, k} $, $ b \in S_{h}^{m',
k'} $, with asymptotic globally homogeneous expansions
\begin{eqnarray*}
a & \sim & \sum_{j\geq 0} a_{-j}, \qquad  a_{-j} \in S_{h}^{m, k +
  \frac{j}{h}}, \ \text{\rm g. h. of degree}\ m-\frac{h}{h+1}k - \frac{j}{h+1}  \\[-5pt]
b & \sim & \sum_{i\geq 0} b_{-i}, \qquad  b_{-i} \in S_{h}^{m', k' +
  \frac{i}{h}}, \ \text{\rm g. h. of degree}\ m'-\frac{h}{h+1}k' - \frac{i}{h+1} .
\end{eqnarray*}
Then $ a \circ b $ is an operator in $ OPS_{h}^{m+m', k+k'} $ with
\begin{multline}
\label{eq:SS} \sigma(a \circ b) - \sum_{s =  0}^{N-1} \sum_{(h+1)|\alpha| + i + j = s} \frac{1}{\alpha!} \sigma\left(\partial_{\xi'}^{\alpha}a_{-j}(x_1, x',
D_{x_1}, \xi') \circ_{x_1} D_{x'}^{\alpha} b_{-i}(x_1, x', D_{x_1}, \xi') \right)
\\
\in S_{h}^{m+m'-N, k+k'}.
\end{multline}
Here $ \circ_{x_1} $ denotes the composition w.r.t. the $x_1$-variable and $\partial_{\xi'}^{\alpha}a_{-j}(x_1, x', D_{x_1}, \xi')$ denotes the
pseudodifferential operator with symbol $\partial_{\xi'}^{\alpha}a_{-j}(x_1, x', \xi_1, \xi')$ quantized in the $(x_1,\xi_1)$ variables.
\end{lemma}
\begin{lemma} [\cite{boutet}, Section 5 and \cite{h2}, Sections 2.2, 2.3]
\label{lemma:H} Let $ a \in H_{h}^{m} $, $ b \in H_{h}^{m'} $ and
$ \lambda \in S_{1,0}^{m''}(\R^{n-1}_{x'}\times\R^{n-1}_{\xi'}) $
with homogeneous asymptotic expansions
\begin{eqnarray*}
a & \sim & \sum_{j\geq 0} a_{-j/(h+1)}, \qquad  a_{-j/(h+1)} \in H_{h}^{m -
  \frac{j}{h+1}}, \ \text{\rm g. h. of degree}\ m- \frac{j}{h+1}
 \\[-5pt]
b & \sim & \sum_{i\geq 0} b_{-i/(h+1)}, \qquad  b_{-i/(h+1)} \in H_{h}^{m' -
  \frac{i}{h+1}},  \ \text{\rm g. h. of degree}\ m' - \frac{i}{h+1} \\[-5pt]
\lambda & \sim & \sum_{\ell\geq 0} \lambda_{-\ell/(h+1)}, \qquad \lambda_{-\ell/(h+1)} \in S_{1, 0}^{m''- \frac{\ell}{h+1}},  \ \text{\rm homogeneous of degree}\
m''- \frac{\ell}{h+1}
\end{eqnarray*}
Then
\begin{itemize}
\item[(i)]{} $ a \circ b^{*} $ is an operator in $
OP\mathcal{H}_{h}^{m+m'- \frac{1}{h+1}}$ with
\begin{multline}
\label{eq:hh*} \sigma(a \circ b^{*})(x_1, x', \xi_1, \xi') -
e^{-ix_1 \xi_1} \sum_{s =
  0}^{N-1} \sum_{(h+1)|\alpha| + i + j = s}
\frac{1}{\alpha!} \partial_{\xi'}^{\alpha} a_{-j/(h+1)}(x', \xi', x_1) D_{x'}^{\alpha} \hat{\bar{b}}_{-i/(h+1)}(x', \xi', \xi_1)
\\
\in \mathcal{H}_{h}^{m + m' - \frac{1}{h+1} - \frac{N}{h+1}},
\end{multline}
where the Fourier transform in $D_{x'}^{\alpha} \hat{\bar{b}}_{-i/(h+1)}$ is taken w.r.t. the $ x_1 $-variable.
\item[(ii)]{} $ b^{*} \circ a $ is an operator in $
OPS_{1,0}^{m+m'- \frac{1}{q}}(\R^{n-1}_{x'}) $ with
\begin{multline}
\label{eq:h*h} \sigma(b^{*} \circ a)(x',\xi') -  \sum_{s =
0}^{N-1} \sum_{(h+1)|\alpha| + j
  + i = s}
\frac{1}{\alpha!} \int \partial_{\xi'}^{\alpha} \bar{b}_{-i/(h+1)}(x',\xi',x_1) D_{x'}^{\alpha} a_{-j/(h+1)}(x',\xi',x_1) dx_1
\\
\in S_{1,0}^{m + m' - \frac{1}{h+1} -
\frac{N}{h+1}}(\R^{2(n-1)}_{(x',\xi')}).
\end{multline}
\item[(iii)]{} $ a \circ \lambda $ is an operator in $ OPH_{q}^{m
+ m''} $. Furthermore its asymptotic expansion is given by
\begin{equation}
\label{eq:hs} \sigma(a \circ \lambda) - \sum_{s=0}^{N-1}
\sum_{(h+1)|\alpha|  + j + \ell =
  s} \frac{1}{\alpha!}
\partial_{\xi'}^{\alpha}a_{-j/(h+1)}(x', \xi',
x_1) D_{x'}^{\alpha}\lambda_{-\ell/(h+1)}(x', \xi') \in H_{h}^{m + m'' - \frac{N}{h+1}}.
\end{equation}
\end{itemize}
\end{lemma}
\begin{lemma}
\label{lemma:SH} Let $ a(x, D) $ be an operator in the class $
OPS^{m,k}_{h} $ and $ b(x', \xi', D_{x_1}) \in OPH_{h}^{m'} $ with
g.h. asymptotic expansions
\begin{eqnarray*}
a & \sim & \sum_{j\geq 0} a_{-j}, \qquad  a_{-j} \in S_{h}^{m, k +
  \frac{j}{h}}, \ \text{\rm g. h. of degree}\ m-\frac{h}{h+1}k - \frac{j}{h+1}  \\[-5pt]
b & \sim & \sum_{i\geq 0} b_{-i/(h+1)}, \qquad  b_{-i/(h+1)} \in H_{h}^{m'-
  \frac{i}{h}}, \ \text{\rm g. h. of degree}\ m'- \frac{i}{h+1} .
\end{eqnarray*}
Then $ a \circ b \in OPH_{q}^{m + m' - k \frac{q-1}{q}}$ and has a
g.h. asymptotic expansion of the form
\begin{multline}
\label{eq:SH} \sigma(a \circ b) -  \sum_{s = 0}^{N-1}
\sum_{(h+1)|\a| + i + j = s} \frac{1}{\a!}
\partial_{\xi'}^{\a}a_{-j}(x_1, x', D_{x_1}, \xi')
\big(D_{x'}^{\a}b_{-i/(h+1)}(x',\xi',\cdot)\big)
\\
\in H_{h}^{m + m' - k \frac{h}{h+1} - \frac{N}{h+1} }.
\end{multline}
\end{lemma}
\begin{lemma}
\label{lemma:HS} Let $ a(x, D) $ be an operator in the class $
OPS^{m,k}_{h}$, $ b^{*}(x',D_{x'},x_1) \in {OPH_{h}^{*}}^{m'} $
and $ \lambda(x', D_{x'}) \in OPS_{1,0}^{m''}(\R^{n-1}_{x'}) $
with homogeneous asymptotic expansions
\begin{eqnarray*}
a & \sim & \sum_{j\geq 0} a_{-j}, \qquad  a_{-j} \in S_{h}^{m, k +
  \frac{j}{h}}, \ \text{\rm g. h. of degree}\ m-\frac{h}{h+1}k - \frac{j}{h+1}  \\[-5pt]
b & \sim & \sum_{i\geq 0} b_{-i/(h+1)}, \qquad  b_{-i/(h+1)} \in H_{h}^{m'-
  \frac{i}{h}}, \ \text{\rm g. h. of degree}\ m'- \frac{i}{h+1} \\[-5pt]
\lambda & \sim & \sum_{\ell\geq 0} \lambda_{-\ell/(h+1)}, \qquad \lambda_{-\ell/(h+1)} \in S_{1, 0}^{m''- \frac{\ell}{h+1}},  \ \text{\rm homogeneous of degree}\
m''- \frac{\ell}{h+1}
\end{eqnarray*}
Then
\begin{itemize}
\item[(i)]{} $ b^{*}(x',D_{x'},x_1) \circ a(x, D) \in
{OPH^{*}_{h}}^{m + m' - \frac{h}{h+1} k}  $ with g.h. asymptotic
expansion
\begin{multline}
\label{eq:SH*} \sigma(b^{*} \circ a) -  \sum_{s = 0}^{N-1} \sum_{(h+1)|\a| + i + j = s} \frac{1}{\a!} \left(D_{x'}^{\a}\overline{a_{-j}}(x_1, x',
  D_{x_1}, \xi')\right)^{*}
(\partial_{\xi'}^{\a}\overline{b_{-i/(h+1)}}(x',\xi',\cdot))
\\
\in {H_{h}}^{m + m' - k \frac{h}{h+1} - \frac{N}{h+1} }.
\end{multline}
\item[(ii)]{} $ \lambda(x', D_{x'}) \circ b^{*}(x',D_{x'},x_1)
\in {OPH_{h}^{*}}^{m' + m''} $ with g.h. asymptotic expansion
\begin{equation}
\label{eq:S10H*} \sigma(\lambda \circ b^{*}) - \sum_{s=0}^{N-1}
\sum_{(h+1)|\alpha| + i +
  \ell = s}
\frac{1}{\alpha!} \partial_{\xi'}^{\alpha}\lambda_{-\ell/(h+1)}(x', \xi') D_{x'}^{\alpha} \overline{b_{-i/(h+1)}}(x',\xi', x_1) \in H_{h}^{m' + m'' -
\frac{N}{h+1}}.
\end{equation}
\end{itemize}
\end{lemma}
The proofs of Lemmas \ref{lemma:SS} --\ref{lemma:SH} are obtained
with a $ h $-variation of the calculus developed by Boutet de
Monvel and Helffer, \cite{boutet}, \cite{h2}. The proof of Lemma
\ref{lemma:HS} is performed by taking the adjoint
$$b^{*}(x',D_{x'},x_1) \circ a(x,D)=\Big(a(x,D)^\ast\circ b^{*}(x',D_{x'},x_1)^\ast\Big)^\ast$$
and using Lemma \ref{lemma:phi*} and \ref{lemma:SH}.\\
We complete this appendix by showing the continuity properties of
the operators defined above.
\begin{lemma}\label{Sobolev}\hspace{50cm}
\begin{itemize}
\item[$(a)$ ]{} Let $ a(x,\xi) \in S_{h}^{m,k} $, properly
supported, with $ k \leq 0 $. Then $ a(x,D) $ is continuous from $
H^{s}_{loc}(\R^{n}) $ to $ H_{loc}^{s - m + k
\frac{h}{h+1}}(\R^{n}) $. \item[$(b)$ ]{} Let $ \phi(x',\xi', x_1)
\in H^{m+ \frac{1}{2(h+1)}}_{h} $, properly supported. Then $
\phi(x',D_{x'}, x_1) $ is continuous from $ H^{s}_{loc}(\R^{n-1})
$ to $ H_{loc}^{s - m}(\R^{n}) $. Moreover $ \phi^{*}(x',D_{x'},
x_1) $ is continuous from $ H^{s}_{loc}(\R^{n}) $ to $ H_{loc}^{s
- m}(\R^{n-1}) $.
\end{itemize}
\end{lemma}

\end{document}